\documentclass[a4paper,10pt,oneside,reqno]{amsart}

\usepackage{amssymb,latexsym,amsthm}

\theoremstyle{plain}
\newtheorem{thm}{Theorem}[section]
\newtheorem*{b2}{Bonatti's Theorem}
\newtheorem{lemma}[thm]{Lemma}
\newtheorem*{ml}{Main Lemma}

\numberwithin{equation}{thm}

\theoremstyle{definition}
\newtheorem{defn}{Definition}[section]

\newcommand{\sig}{\,$\Sigma$\, }

\begin{document}

\title[Commuting diffeomorphisms]{A note on commuting\\ diffeomorphisms  on
surfaces}

\author{S. Firmo}



\date{March 6, 2005}

\address{Sebasti\~ao Firmo\\ Universidade Federal Fluminense\\Instituto de Matem\'a\-tica\\
Rua M\'ario Santos Braga \\ 24020\,-140 Niter\'oi, RJ\\Brasil
\vglue10pt}

\email{firmo@mat.uff.br}

\maketitle


\begin{abstract} Let \,$\Sigma$\, be a closed surface with nonzero Euler
characteristic.  We prove the existence of an open neighborhood
\,$\mathcal{V}$\, of the identity map of  \,$\Sigma$\, in the \,$C^1$-topology 
with the following property: if
\,$G$\, is an abelian subgroup of \,$\textrm{Diff}^{\hskip1pt 1}(\Sigma)$\, generated by any family of
elements in \,$\mathcal{V}$\, then the elements of \,$G$\, have common fixed points. 
This result generalizes a similar result due to Bonatti and  announced in his paper
\emph{Diff\'eomorphismes  commutants des surfaces et stabilit\'e des fibrations en tores}.
\end{abstract}


\medskip
\medskip

\thispagestyle{empty}


\section{Introduction}
\vskip3mm


Bonatti has proven in \cite{b2} the following result.

\medskip

\begin{b2} 

Let \,$\Sigma$ be a closed surface with nonzero Euler characteristic. 
Fixed \,$k\in\mathbb{Z}^+$, there is an open \,$C^1$-neighborhood \,$\mathcal{U}_k$\,
of the identity map of 
\,$\Sigma$ satisfying the following\,$:$
if
\,$G$\, is an abelian subgroup of \,$\textrm{\rm Diff}^{\hskip1pt 1}(\Sigma)$\, 
generated by \,$k$\,  elements in 
\,$\mathcal{U}_{k}$\, then for some \,$p\in\Sigma$\, we have \,$f(p)=p$\, for all \,$f\in G$.

\end{b2}

\medskip

We remark that in the above theorem the size of the neighborhood \,$\mathcal{U}_k$\, depends on the
number \,$k$\, of generators
of the abelian group \,$G$\, unless \,$\Sigma$\, is the \,$2$-sphere \,$S^2$ or the
projective plane \,$\mathbb{RP}^2$. 
In fact, the  cases \,$S^2$\, and \,$\mathbb{RP}^2$\, were  treated
by Bonatti in
\cite{b1} and in these cases  the neighborhood of the identity map can be chosen to be uniform.

The purpose of this paper is to prove that even when  
\,$\Sigma$\, is different from \,$S^2$\, and \,$\mathbb{RP}^2$,  there  exists a 
distinguished  \,$C^1$-neighborhood 
of the identity map of \,$\Sigma$\,  where the above theorem
holds regardless of the number of generators of the group \,$G$. Precisely, we prove the
following theorem.

\medskip

\begin{thm}\label{TC}  
 
Let \,$\Sigma$\, be a closed  surface with nonzero
Euler characteristic. Then there exists an open \,$C^1$-neighborhood
\,$\mathcal{V}$   of the identity map of \,$\Sigma$  having the
following property\,$:$ if 
\,$G$\, is an abelian subgroup of \,$\textrm{\rm Diff}^{\hskip1pt 1}(\Sigma)$\,
 generated by any family of elements in 
\,$\mathcal{V}$\, then 
for some \,$p\in\Sigma$\, we have \,$f(p)=p$\, for all \,$f\in G$.

\end{thm}

\medskip

With regard to the techniques of Bonatti's paper \cite{b2} we observe that he provides 
a neighborhood \,$\mathcal{U}_2$\, of the
identity map of 
\,$\Sigma$\, such that two commuting diffeomorphisms in this neighborhood have common fixed points.
To guarantee the existence of a common fixed point for three commuting diffeomorphisms,  he needs  to
shrink the neighborhood \,$\mathcal{U}_2$\, to a certain neighborhood \,$\mathcal{U}_3$. Similarly, 
\,$\mathcal{U}_3$\, have to be shrink if there is four or more diffeomorphisms.
Our argument consists of showing that the above mentioned sequence 
\,$\mathcal{U}_2,\mathcal{U}_3,\ldots$\, of neighborhoods actually stabilizes at some integer
depending only on the topology of the surface.

Since we are dealing  with abelian groups generated by diffeomorphisms  close to the
identity map their lifts to a double covering space still form an abelian group with generators
close to the identity.  
Therefore, by using  the double covering of the orientations of \sig we conclude that to prove
Theorem \ref{TC} it suffices to prove it  for orientable closed surfaces
with nonzero Euler characteristic.

We close the introduction with the following question.

Does our theorem hold for  homeomorphisms \,$C^0$\, close to the identity map of a closed surface 
\,$\Sigma$? It seems to the author that this question remains open even for 
\,$C^1$-diffeomorphisms that are \,$C^0$ close to the identity map.
Handel, in \cite{ha}, proves the existence of common fixed points for two commuting
homeomorphisms of \,$S^2$\, which are, for example, \,$C^0$ close to the identity map. Moreover,
as to higher genus  surfaces,    he proves that two orientation preserving \,$C^1$-diffeomorphisms
\,$f$\, and \,$g$ which commute have at least as many common fixed points as \,$F$\, and \,$G$\, do,
provided that \,$F$\, and \,$G$\, are  pseudo-Anosov. Here,
\,$F$\, and \,$G$\, stand for  the homeomorphisms respectively  obtained from \,$f$\, and \,$g$\,
by the Thurston Classification Theorem for surface homeomorphisms.


\vskip35pt

\section{Notations and definitions}\label{sec:2} 
\vskip3mm


From now on  \,$\Sigma$\, will be a closed connected oriented surface embedded in the Euclidean
space
\,$\mathbb{R}^3$\, endowed with the usual norm denoted by \,$\| \cdot \|$. The distance
in \,$\Sigma$\, associated to the induced Riemannian metric from \,$\mathbb{R}^3$\, will be
denoted by \,$d$\, and we denote by \,${\rm B}(x\,;\rho)$\, the closed $2$-ball centered at \,$x$\,
with radius
\,$\rho \in(0\,,\infty)$\, with respect to  \,$d$.

The \,$C^1$-norm of a given \,$C^1$-map
 \,$\varphi:\Sigma\rightarrow \mathbb{R}^3$\, is defined 
by 
$$\|\varphi\|_1=\sup_{x\in \Sigma}\Big\{\|\varphi(x)\|+
\sup_{v\in T_{x}\Sigma\,;\,\|v\|=1}\|D\varphi(x)\cdot v\|\Big\}.$$

\noindent
Defined on  the real vector space consisting of \,$C^1$-maps from \,$\Sigma$\, to \,$\mathbb{R}^3$,
this norm induces
a distance between \,$C^1$-maps \,$\varphi,\psi:\Sigma\rightarrow\mathbb{R}^3$ 
given by \,$\|\varphi-\psi\|_1$.
The distance between 
two \,$C^1$-diffeomorphisms \,$f,h$\, of \,$\Sigma$\, 
is, by definition, the distance between \,$f$\, and \,$h$\, as \,$C^1$-maps from 
\,$\Sigma$\, to \,$\mathbb{R}^3$. This distance on the space of all \,$C^1$-diffeomorphisms of 
\,$\Sigma$\, defines the so-called \,$C^1$-topology and
the group of \,$C^1$-diffeomorphisms of \,$\Sigma$\, endowed with this topology  will be
denoted by 
\,$\textrm{\rm Diff}^{\hskip1pt 1}(\Sigma)$. 

Given a list  \,$\mathcal{H}$\, 
of elements in 
\,$\textrm{\rm Diff}^{\hskip1pt 1}(\Sigma)$\,, let \,${\rm Fix}(\mathcal{H})$\, 
denote the set of 
common fixed points  of the elements of  \,$\mathcal{H}$. In other words, 
\,${\rm Fix}(\mathcal{H})=\bigcap_{f\in\mathcal{H}} {\rm Fix}(f)$\, 
where \,${\rm Fix}(f)$\, is the set of fixed
points of \,$f$.

The \,\emph{positive semi-orbit}\, of  \,$p\in \Sigma$\, by a  diffeomorphism
\,$f$\,  is the set \,$\mathcal{O}_p^+(f)=\{f^n(p) \ ;
\ n\geq 0\}$. Its closure in \,$\Sigma$\, will be denoted by
\,$\overline{\mathcal{O}_p^+(f)}$\,. We  say that \,$p\in\Sigma$\, is a
\,\emph{$\omega$-recurrent} point  for \,$f$\, if
\,$p$\, is the limit of some subsequence of \,$\big(f^n(p)\big)_{n\geq0}$\,.

Let \,$a,b\in\Sigma$\, be such that \,$d(a,b)$\, is smaller than the injectivity radius of the
exponential map associated to the metric of \,$\Sigma$. For such \,$a\,,b$\, 
we denote by
\,$[\,a\,,b\,]$\, the  oriented geodesic arc  joining \,$a$\, to \,$b$\, which is contained in
the disc of injectivity
centered at \,$a$.

Fix \,$\kappa>0$\, such that \,$d(a,b)$\, is smaller than the injectivity radius of the
exponential map whenever \,$\|b-a\|<\kappa$.

Now, given  \,$f\in\textrm{Diff}^{\hskip1pt 1}(\Sigma)$\, such that
\,$\|f-\text{Id}\|_1<\kappa$\,, let \,$X_f$\, be the standard vector field on
\sig{} associated to \,$f$\, as follows: 
\,$X_f(x)$\, is the tangent vector to the geodesic segment \,$[\,x\,,f(x)\,]$\, at the point \,$x$\,
with orientation given by  the  orientation of \,$[\,x\,,f(x)\,]$\,  whose norm is equal to
\,$\|f(x)-x\|$. The singular set of
\,$X_f$\, is the set 
\,${\rm Fix}(f)$.

Moreover, let  \,$p\in\Sigma-{\rm Fix}(f)$\,. Following Bonatti \cite{b2},  a 
piecewise geodesic simple closed curve
\,$\Gamma^{\hskip1pt p}_{\!\!f}$\, on \sig{} is said to be
\emph{supported by}
\,$\mathcal{O}_p^+(f)$\, if each geodesic arc of \,$\Gamma^{\hskip1pt p}_{\!\!f}$\, is
contained in some segment \,$[f^i(p)\,,f^j(p)]$\,   satisfying 
$$d\big(f^i(p)\,,f^j(p)\big)\leq \mbox{$\frac{3}{2}$} \, d\big(f^i(p)\,,f^{i+1}(p)\big) \ \ 
\mbox{where} \ \ i\,,j\geq0.$$

Consider \,$0<\alpha<\pi$. We say that \,$\Gamma^{\hskip1pt p}_{\!\!f}$\, is  
  \,$\alpha$\emph{-tangent} to the vector field \,$Y$ provided that the following conditions
are satisfied:  \,$Y$ has no singularities along \,$\Gamma^{\hskip1pt p}_{\!\!f}$\, and
for one of the two possible orientations of \,$\Gamma^{\hskip1pt p}_{\!\!f}$\, and for 
each point \,$x$\, of a  geodesic arc \,$\Omega$\, of \,$\Gamma^{\hskip1pt p}_{\!\!f}$\, the angle
 between the vector 
\,$Y(x)$\, and the unitary tangent vector to \,$\Omega$\, at \,$x$\, induced by the
orientation of  \,$\Gamma^{\hskip1pt p}_{\!\!f}$\, is less than \,$\alpha$\, for all
\,$x\in\Omega$.

Furthermore, let \,$\gamma^{\,p}_{f}$\, denote the  curve obtained by concatenating the 
segments
\,$[f^i(p)\,,f^{i+1}(p)]$\, for \,$i\geq0$.

The following topological result about  compact surfaces will be very important in our proofs. 
It will be used in the proof of Theorem \ref{TPF} and in the proof of the Main Lemma.

Let \,$\Sigma\subset\mathbb{R}^3$\, be an oriented connected closed  surface.
Then  there is an integer
\,$N\geq 4$\, with the following property:

Given any compact connected surface \,$\mathcal{S}\subset\Sigma$\,  with boundary 
 such that each connected component of its boundary  is not null homotopic in \,$\Sigma$, one has:

\begin{itemize}

\item 
The number of connected components of  \,$\partial\mathcal{S}$\, is less than \,$N$;

\item 
If \,$\alpha_1,\ldots,\alpha_{N}\subset{\rm Int}(\mathcal{S})$\, is a list of pairwise disjoint
simple closed curves  then there are two distinct curves
\,$\alpha_i\,,\alpha_j$\,  in that list which are homotopic in \,$\mathcal{S}$.
Consequently,
\begin{itemize}
\item [--]
either each one of these two curves bounds  disks embedded in \,$\mathcal{S}$\,;

\item[--]
or   these two curves bound a cylinder embedded in \,$\mathcal{S}$.

\end{itemize}

\end{itemize}

We remark that the constant \,$N\geq 4$\,  considered in the last paragraph  will 
frequently be used in all this paper.


\vskip35pt

\section{Some known results}\label{sec:3}
\vskip3mm


In this section we recall some technical results from
\cite{b2}  which play a key role in this note.

Given \,$\epsilon>0$\, let 
\,$\{\mathcal{V}_k(\epsilon)\}_{k\geq1}$\, be a decreasing nested sequence 
of open  neighborhoods of the identity map of
\,$\Sigma$\, in the  \,$C^1$-topology inductively defined as follows:

\begin{itemize}

\item 
$\mathcal{V}_1(\epsilon) =\big\{f \in\mbox{\rm Diff}^{\hskip1pt 1}(\Sigma) \ ; \
\|f-\text{Id}\|_1<\epsilon  \big\}  \nonumber$;

\item
Fixed \,$\mathcal{V}_k(\epsilon)$\, for some positive integer \,$k$\, we choose
\,$\mathcal{V}_{k+1}(\epsilon)$\,
so that the following holds: 
$$\text{if} \ \
f_1,\ldots,f_{2N+1}\in\mathcal{V}_{k+1}(\epsilon)
\ \ \text{then} \ \
f_1 \text{\tiny$\circ$} f_2 \, \text{\tiny$\circ$} \cdots
\text{\tiny$\circ$} f_{2N} \text{\tiny$\circ$}
f_{2N+1}\in\mathcal{V}_{k}(\epsilon)$$ where 
\,$f_1 \text{\tiny$\circ$} f_2 \, \text{\tiny$\circ$}\cdots 
\text{\tiny$\circ$} f_{2N} \text{\tiny$\circ$} f_{2N+1}$\, stands for the
composition of  maps.

\end{itemize}

We always assume that \,$\epsilon<\kappa$\, so that we can guarantee that 
\,$[\,p\,,f(p)]$\, is well defined whenever \,$f\in\mathcal{V}_{1}(\epsilon)$.


\medskip

\begin{lemma}\label{bona2}{\rm\textbf{(Bonatti)}}
There exists  \,$0<\epsilon_1<\kappa$\, such that
every pair of elements  \,$f,h\in\mbox{\rm Diff}^{\hskip1pt 1}(\Sigma)$\,
satisfies the following$\,:$
\begin{enumerate}

\item\label{bona1}
If  \,$f\in\mathcal{V}_1(\epsilon_1)$\, and 
 \,$f(p)\neq p$\, then \,$f$\, does not have fixed points in the
ball  \,$\mbox{\rm B}\big(p\,;4\,d(p\,,f(p))\big)$. In particular, \,$f$\, 
does not have fixed points along curves supported by \,$\mathcal{O}_p^+(f)$\, if such curves
exist.

\item\label{char}
If \,$f\in\mathcal{V}_1(\epsilon_1)$\, and   \,$p\in \Sigma-{\rm Fix}(f)$\, is a
\,$\omega$-recurrent point of
\,$f$\, then there
is a  curve \,$\Gamma^{\hskip1pt p}_{\!\!f}$\, supported by \,$\mathcal{O}_p^+(f)$\, and
\,$\frac{\pi}{10}$-tangent to \,$X_f$.

\item\label{fixp}
If \,$f,h\in\mathcal{V}_{2}(\epsilon_1)$\,  commute  and  
\,$p_i\in{\rm Fix}(h \, \mbox{\tiny$\circ$} \, f^i)-{\rm Fix}(f)$\, then
\,$h \, \mbox{\tiny$\circ$} \, f^j$\,  has no fixed points in the ball  \,${\rm
B}\big(p_i\,;4\,d(p_i\, ,f(p_i)\big)$\, where \,$j\neq i$\, and
\,$i\,,j\in\{0,\ldots,2N\}$. In particular, \,$h \, \mbox{\tiny$\circ$} \, f^j$\, does not have fixed
points along curves supported by \,$\mathcal{O}^{+}_{p_i}(f)$.

\item\label{fixp2}
If \,$f,h\in\mathcal{V}_{2}(\epsilon_1)$\,  commute,
\,$0<d(p_i\,,p)\leq\frac{3}{2}\,d(p_i\,,f(p_i))$\, and 
 the segment \,$[\,p_i\,,p\,]$\, 
is  \,$\frac{\pi}{10}$-tangent to the vector field \,$X_f$\, then  
the segment \,$[\,p_i\,,p\,]$\, is \,$\frac{2\pi}{5}$-tangent
to the vector field 
\,$X_{h \, \mbox{\tiny$\circ$} \, f^j}$\, for all \,$j\neq i$. Here 
\,$i\,,j\in\{0,1,\ldots,2N\}$,
\,$p\in\Sigma$\, and 
 \,$p_i\in{\rm Fix}(h \, \mbox{\tiny$\circ$} \, f^i)-{\rm Fix}(f)$.

\end{enumerate}
\end{lemma}


The simple closed curve \,$\Gamma^{\hskip1pt p}_{\!\!f}$\, obtained in item
(\ref{char}) of Lemma \ref{bona2} will be called  \,\emph{character curve}\, of \,$f$\, at \,$p$.

The results listed in Lemma \ref{bona2} are proven in \cite{b2}  for the
integers \,$i\,,j\in\{0,\ldots,N\}$\, where \,$N$\, was defined at the end of the last section.
Nonetheless, it is easy to see that the proofs in \cite{b2} also work for an arbitrary  
positive integer modulo reducing \,$\epsilon_1$.


\vskip35pt

\section{Preparing the proof of Theorem \ref{TC}}
\vskip3mm


In this section we shall establish three lemmas which will be used in the next section to prove 
Theorem \ref{TC}.
The proof of the first lemma
is, in fact,  contained in the proof of Lemma 4.1 of \cite[pages~67--68]{b1}. We repeat 
the arguments 
here because our hypothesis are not exactly the same  as those used in Bonatti's paper.
Besides, in section \ref{prooftechlem}, the argument below will be further 
adapted to apply to more general and technical situations. 
Lemma \ref{fixp4A} below says that  curves supported by positive
semi-orbits are disjoint under appropriate conditions.

From now on the neighborhood \,$\mathcal{V}_k(\epsilon_1)$\, will be  denoted simply by 
\,$\mathcal{V}_k$\, for all \,$k\in\mathbb{Z}^+$\, where \,$\epsilon_1$\, 
is always given by Lemma \ref{bona2}.

\medskip

\begin{lemma}\label{fixp4A} 
Let \,$f,h\in\mathcal{V}_{1}$\, be commuting diffeomorphisms 
and let 
\,$\Gamma^{\hskip1pt p}_{\!\!f},\Gamma^{\hskip1pt q}_{\!\!h}$\, be curves supported
respectively by 
\,$\mathcal{O}^{+}_{p}(f)\,,\mathcal{O}^{+}_{q}(h)$\,   where 
$$p\in {\rm Fix}(h)-{\rm Fix}(f) \quad  \text{and} \quad  q\in {\rm Fix}(f)-{\rm Fix}(h).$$
Then we have$\,:$

\begin{itemize}

\item  
$\gamma^{\hskip1pt p}_{f} \, \cap \, \Gamma^{\hskip1pt q}_{\!\!h}=\emptyset \, ;$

\item  \,$\Gamma^{\hskip1pt {p}}_{\!\!f}\cap\Gamma^{\hskip1pt q}_{\!\!h}=\emptyset \,;$

\item
$\overline{\mathcal{O}^{+}_{p}(f)} \, \cap \, 
\Gamma^{\hskip1pt q}_{\!\!h}=\emptyset$.

\end{itemize}

Moreover,
if \,$\rho>0$\,  is such that
$$d\big(x\,,f(x)\big)\,,d\big(y\,,h(y)\big) \geq \rho \ \ , \ \ \forall \, x\in 
\mathcal{O}^{+}_{p}(f) \ \ \text{and} \ \ \forall \, y\in \mathcal{O}^{+}_{q}(h)  $$
then  \,$d(\Gamma^{\hskip1pt p}_{\!\!f},\Gamma^{\hskip1pt q}_{\!\!h})\geq \rho$.

\end{lemma}


\begin{proof}

Let us first  prove  the second item. To do this we suppose for a contradiction that
\,$\Gamma^{\hskip1pt {p}}_{\!\!f}\cap \, \Gamma^{\hskip1pt q}_{\!\!h}\neq\emptyset$.
Thus, there exist integers \,$m,n,k,l\geq0$\, such that
$$[f^m(p)\,,f^n(p)]\cap [h^k(q)\,,h^l(q)]\neq\emptyset$$
where 
\begin{equation}\label{sec4:lem4.1:sist1}
\begin{split}
d\big(f^m(p)\,,f^n(p)\big) & \leq \mbox{$\frac{3}{2}$} \, d\big(f^m(p)\,,f^{m+1}(p)\big)\\
d\big(h^k(q)\,,h^{l}(q)\big) & \leq \mbox{$\frac{3}{2}$} \, d\big(h^k(q)\,,h^{k+1}(q)\big).
\end{split}
\end{equation}

\noindent
By the triangle inequality we obtain:  
\begin{align}
 d\big(f^m(p)\,,h^k(q)\big) & \leq   
d\big(f^m(p)\,,f^n(p)\big)+d\big(h^k(q)\,,h^{l}(q)\big) \nonumber \\
 & \leq   
\mbox{$\frac{3}{2}$} \, d\big(f^m(p)\,,f^{m+1}(p)\big)+\mbox{$\frac{3}{2}$} \,
d\big(h^k(q)\,,h^{k+1}(q)\big) 
\nonumber\\
& \leq 3\max\Big\{d\big(f^m(p)\,,f^{m+1}(p)\big) \,,\,
d\big(h^k(q)\,,h^{k+1}(q)\big)\Big\}. 
\nonumber
\end{align}
Therefore, we have the following two possibilities:
\begin{itemize}
\item [--]
either \,$h^k(q)$\, is  in the   ball 
 \,$\textrm{B}\big(f^m(p)\,;3\,d (f^m(p)\,,f^{m+1}(p))\big)$\, which is impossible
by item (\ref{bona1}) of Lemma \ref{bona2}  since        
the map \,$f$\, has no fixed points in 
 \,$\textrm{B}\big(f^m(p)\,;3\,d (f^m(p)\,,f^{m+1}(p))\big)$. Note that \,$h^k(q)\in{\rm Fix}(f)$\, 
which follows from the commutativity;

\item[--]
or \,$f^m(p)$\, is in the ball
\,$\textrm{B}\big(h^k(q)\,;3\,d(h^k(q)\,,h^{k+1}(q))\big)$\,   
which is impossible by the same reason.
\end{itemize}

\noindent
This finish the  proof of the second item.

The reader will notice that 
the above arguments  prove also the first
item. 
The  last item  follows from observing that 
\,$\overline{\mathcal{O}^{+}_{p}(f)}\subset\text{Fix}(h)$\, and \,$h$\, is free of fixed points over
\,$\Gamma^{\hskip1pt q}_{\!\!h}$\, thanks to item (\ref{bona1}) of Lemma \ref{bona2}.

To prove the second part of the lemma let us suppose for a contradiction that  
\,$d(\Gamma^{\hskip1pt p}_{\!\!f},\Gamma^{\hskip1pt q}_{\!\!h})< \rho$. Then
there exist \,$m,n,k,l\geq0$\,  satisfying (\ref{sec4:lem4.1:sist1})  and two
points 
$$a\in[f^m(p)\,,f^n(p)] \ \ \text{and} \ \ b\in[h^k(q)\,,h^l(q)]$$
such that \,$d(a\,,b)<\rho$. Therefore,
\begin{align}
 d\big(f^m(p)\,,h^k(q)\big) & \leq   
d\big(f^m(p)\,,f^n(p)\big)+d(a\,,b)+d\big(h^k(q)\,,h^{l}(q)\big) \nonumber \\
 & \leq    
\mbox{$\frac{3}{2}$} \, d\big(f^m(p)\,,f^{m+1}(p)\big)+\rho+\mbox{$\frac{3}{2}$} \,
d\big(h^k(q)\,,h^{k+1}(q)\big) 
\nonumber\\
& \leq  4\max\Big\{d\big(f^m(p)\,,f^{m+1}(p)\big) \,,\,
d\big(h^k(q)\,,h^{k+1}(q)\big)\Big\}. 
\nonumber
\end{align}
Now, we finish the proof by using exactly the same arguments  used to prove that
\,$\Gamma^{\hskip1pt {p}}_{\!\!f}\cap\Gamma^{\hskip1pt q}_{\!h}=\emptyset$.
\end{proof}

\medskip

The next lemma is a version of Lemma 5.1 in \cite[page~69]{b1} for the surface \,$\Sigma$.
It is a version of Bonatti's Theorem in \cite{b2} for a special kind of boundary.

Here we use the notation 
\,$x\in\text{Fix}(f_1,\ldots,{\widehat f}_{\!\lambda}, \ldots,f_m)$\, 
to mean that \,$x\in\text{Fix}(f_i)$\, for all \,$i\in\{1,\ldots,m\}$\, and \,$i\neq \lambda$.

\medskip

\begin{lemma}\label{V.n}

Let \,$f_{1},\ldots,f_{3N}\in\mathcal{V}_{3N+1}$\,
be commuting diffeomorphisms and
let    \,$\mathcal{S}\subset\Sigma$\, be a compact connected surface with
\,$\chi(\mathcal{S})\neq0$. Suppose that
any connected component \,$\Gamma$ of 
\,$\partial \mathcal{S}$\, satisfies the following$\,:$
\begin{itemize}
\item  
$\Gamma$ is  not null homotopic 
in \,$\Sigma\,;$ 

\item
$\Gamma$ is a character curve of \,$f_{\lambda}$\, at 
$$p\in {\rm Fix}(f_1 \,, \ldots \,, 
{\widehat f_{\lambda}} \,, \ldots \,, f_{3N}) - {\rm Fix}(f_{\lambda})$$
for some \,$\lambda\in\{1,\ldots,3N\}$.
\end{itemize}

Then    \,$f_1 , \ldots , f_{3N}$\,  
 have common fixed points in \,${\rm Int}(\mathcal{S})$.

\end{lemma}


We notice that this lemma contains the case \,$\mathcal{S}=\Sigma$. 
Moreover, by the definition of the integer \,$N$ we have that the number of connected components of
\,$\partial\mathcal{S}$\, is less than \,$N$.
The proof of this lemma will be deferred to 
section \ref{proofV.n}.

The proof of the next lemma is easily obtained from the proof of Lemma~5.1 in \cite[page~69]{b1}.

\medskip

\begin{lemma}\label{bondiscsup}

Let \,$f_{1},\ldots,f_{n}\in\mathcal{V}_{1}$\,
be commuting diffeomorphisms and let 
 \,$\Gamma^{\hskip1pt {p}}_{\!\!f_n}$\, be a character curve of \,$f_n$\, at
$$p\in {\rm Fix}(f_1 \,, \ldots \,, f_{n-1}) - {\rm Fix}(f_{n}).$$
If \,$\Gamma^{\hskip1pt {p}}_{\!\!f_n}$\, bounds a disc in \,$\Sigma$\, then 
\,$f_1 \,, \ldots \,, f_n$\, have a common fixed point in the interior of that disc.

\end{lemma}



\vskip35pt

\section{Proof of Theorem \ref{TC}}
\vskip3mm


Theorem~\ref{TC} will be obtained as an easy consequence of Theorem~\ref{TPF}
to be stated and proved below. The proof of Theorem~\ref{TPF} is by induction 
and by contradiction. It consists of two parts. In the first part the induction procedure
is initialized (i.e. the first step of the induction is stablished). 
This part will be carried after the proof of  Lemma~\ref{V.n} in section 
\ref{proofV.n}. There, we shall recast
Bonatti's proof \cite{b2} in a more general context (with boundaries) for \,$k$\, 
diffeomorphisms of \,$\Sigma$  where \,$2\leq k\leq 3N$\, (recall that \,$N$\, was defined at the
end of section \ref{sec:3}). This argument will also go by induction and by contradiction.
At each step of the induction procedure over \,$k$, the neighborhood of the identity in question
will be reduced to guarantee the existence of a common fixed point for the diffeomorphisms.
Roughly speaking, the technical reason to reduce the neighborhood is that we need to construct 
\,$N$\, special character curves pairwise disjoint by using \,$k+1$\, diffeomorphisms
\,$f_1,\ldots,f_k,f_{k+1}$\, (once the existence of the common fixed point for \,$k$\, 
diffeomorphisms has been stablished). The construction of these character curves is carried out 
by using the positive semi-orbit by \,$f_{k+1}$\, of common fixed points of the \,$k$\,
diffeomorphisms
\,$f_1,\ldots, f_{k-1}, f_k \mbox{\,\tiny$\circ$}f_{k+1}^i$\, where \,$i\in\{1,\ldots,N\}$.
In view of the diffeomorphism \,$f_k \mbox{\,\tiny$\circ$}f_{k+1}^i$, the neighborhood of the
identity needs to be reduced so as to guarantee that \,$f_k \mbox{\,\tiny$\circ$}f_{k+1}^i$\, 
belongs to the neighborhood obtained in the previous step of the induction (i.e. the case of
\,$k$\, diffeomorphisms). This technical question is already apparent in Bonatti's proof 
of \cite[page~109]{b2}.

In the second part of the proof we have a family  consisting of more than \,$3N$\, commuting
diffeomorphisms. In this case these \,$N$\, special character curves pairwise disjoint will be
obtained through the positive semi-orbit by \,$f_j$\, of the common fixed points 
of the diffeomorphisms 
\,$f_1,\ldots,\widehat{f_j},\ldots, f_{3N+n+1}$\, where \,$f_j$\, 
will be conveniently chosen from the set 
\,$\{f_{N+1},\ldots,f_{3N}\}$.
Finally, with this new construction procedure  of the \,$N$\, pairwise disjoint character curves, we
shall be able  to keep the same neighborhood of the identity for \,$3N$\, diffeomorphisms.

The construction of the character curves carried out in the second part of the proof, which exploits
the existence of a large number of generators, is the essential difference between Bonatti's proof
and the present one. Naturally, in order to apply this strategy, Bonatti's Theorem has to be
extended to a more general settings. Such extension however will be accomplished by using the same
arguments employed in \cite{b1,b2}.

On the other hand, to construct these \,$N$\, character curves supported by appropriate
semi-orbits, it will be necessary to ensure that the semi-orbits remain in 
\,$\text{Int}(\mathcal{S})$. It will also be necessary to guarantee that the corresponding character
curves are disjoint and do not intersect the boundary of \,$\mathcal{S}$.
For all that, Lemma \ref{fixp4A} will be crucial.

As a matter of fact, the strategy used in the second part is implicit in the proof 
of Bonatti's Theorem for \,$S^2$\, in \cite{b1}. In this case, only two diffeomorphisms are needed to
implement the construction of the character curves. As a consequence, the first part of the proof
is superfluous in this case.


\medskip

\begin{thm}\label{TPF}

 Let \,$n\geq0$\, be an integer and let 
$$f_{1},\ldots,f_{N},h_1,\ldots,h_{2N+n}\in\mathcal{V}_{3N+1}$$
 be  commuting   diffeomorphisms. Consider  
a compact connected surface \,$\mathcal{S}\subset\Sigma$\, with \,$\chi(\mathcal{S})\neq0$.
Suppose that  any connected component \,$\Gamma$ of 
\,$\partial \mathcal{S}$\, satisfies the following$\,:$ 

\begin{itemize}
\item 
$\Gamma$ is  not null homotopic 
in \,$\Sigma\,;$ 

\item
$\Gamma$ is a character curve of \,$f_{\lambda}$\, at
$$p\in {\rm Fix}(f_1 \,, \ldots \,, 
{\widehat f_{\lambda}} \,, \ldots \,, f_{N}\,,h_1\,,\ldots\,,h_{2N+n}) - {\rm Fix}(f_{\lambda})$$
for some \,$\lambda\in\{1,\ldots,N\}$.
\end{itemize}

Then  
\,$f_{1},\ldots,f_{N},h_1,\ldots,h_{2N+n}$\,  
 have common fixed points in \,${\rm Int}(\mathcal{S})$.

\end{thm}


\begin{proof}

We argue by induction on  \,$n\geq0$.

The theorem holds for \,$n=0$\, since  this case  reduces to  Lemma~\ref{V.n}.

Now, let us assume that it holds for some integer \,$n\geq0$. Also, let us suppose for a
contradiction that
\,$f_1,\ldots,f_{N},h_1,\ldots,h_{2N+n},h_{2N+n+1}$\, do not   have common fixed
points in 
\,${\rm Int}(\mathcal{S})$.

By assumption, if \,$\partial\mathcal{S}\neq\emptyset$\, then each connected component of 
\,$\partial\mathcal{S}$\, is a character curve  
\,$\Gamma^{p_{\lambda}}_{\!\!f_{\lambda}}$ of \,$f_{\lambda}$\, at
$$p_\lambda\in{\rm Fix}(f_1,\ldots,{\widehat f}_\lambda,\ldots,
f_{N},h_1,\ldots,h_{2N+n},h_{2N+n+1}\big)-{\rm Fix}(f_\lambda)$$
for some \,$\lambda\in\{1,\ldots,N\}$.

For each \,$j\in\{1,\ldots,N\}$\, let us consider the list
$$f_1,\ldots,f_N\,,h_1,\ldots,{\widehat h}_j,\ldots,h_N,\ldots,h_{2N+n+1}$$
of \,$3N+n$\, diffeomorphisms.
From the induction assumption on \,$n$\, we conclude that
there exists a point \,$q_j\in{\rm
Int}(\mathcal{S})$\, such that
\begin{align}\label{}
q_j\in{\rm Fix}(f_1,\ldots,f_N\,,h_1,\ldots,{\widehat h}_j,\ldots,h_{2N+n+1})-{\rm
Fix}(h_j)
\nonumber
\end{align}
since \,$f_1,\ldots,f_N\,,h_1,\ldots,h_{2N+n+1}$\, do not have common
fixed points in  \,$\text{Int}(\mathcal{S})$.
From  Lemma \ref{fixp4A} we know that 
\,$\overline{\mathcal{O}^{+}_{q_j}(h_j)}\subset \text{Int}(\mathcal{S})$. 
On the other hand,  \,$\overline{\mathcal{O}^{+}_{q_j}(h_j)}$\, is invariant by 
\,$h_j$\, and  contained in 
$${\rm Fix}(f_1,\ldots,f_N\,,h_1,\ldots,{\widehat h}_j,\ldots,h_{2N+n+1})-{\rm
Fix}(h_j).$$
Now, thanks to  Zorn's Lemma, we can assume without loss of generality
that \,$q_j$\, is a \,$\omega$-recurrent point for \,$h_j$.

Besides, the maps
\,$f_1,\ldots,f_N\,,h_1,\ldots,h_{2N+n+1}$\, do not have common fixed points over 
\,$\partial\mathcal{S}$\, since  \,$f_{\lambda}$\, has no fixed points 
over \,$\Gamma^{p_{\lambda}}_{\!\!f_{\lambda}}$. Therefore, 
\,$f_1,\ldots,f_N\,,h_1,\ldots,h_{2N+n+1}$\,
 do not have common fixed points in
\,$\mathcal{S}$. Thus, there exists \,$\rho>0$\, satisfying the following
condition:
$$ d\big(x \, ,h_{\ell}(x)\big)\geq\rho$$
for all 
\,$x\in\text{Fix}(f_1,\ldots,f_N\,,h_1,\ldots,{\widehat
   h}_{\ell},\ldots,h_{2N+n+1})\cap\mathcal{S}$\, and for all 
 \,$\ell\in\{1,\ldots,2N\}$.

Let \,$\delta>0$\, be such that the area of any disk of radius \,$\rho/3$\, contained in 
\,$\Sigma$\, is greater than \,$\delta$.

From  Lemma \ref{fixp4A}  we have that the character curves
\begin{align}\label{sec6:N1}
\Gamma^{\hskip1pt q_1}_{\!\!h_1},\ldots,\Gamma^{\hskip1pt q_N}_{\!\!h_N}
\end{align}
are contained in \,$\text{Int}(\mathcal{S})$
and  the distance between any two distinct curves of the above list
is greater than or equal to \,$\rho$.

On the other hand it follows from the topology of \,$\mathcal{S}$\, the existence of 
two distinct curves \,$\Gamma^{\hskip1pt q_i}_{\!\!h_i}$\, and \,$\Gamma^{\hskip1pt
q_j}_{\!\!h_j}$\,  in the list \ref{sec6:N1} which are
homotopic. Furthermore, they cannot  bound any disc in \,$\mathcal{S}$\, since in that case, 
it would  follow from 
Lemma \ref{bondiscsup} that
\,$f_1,\ldots,f_N\,,h_1,\ldots,h_{2N+n+1}$\,  have  common fixed points in the interior of the
disc in question. This is however impossible.

Consequently, 
\,$\Gamma^{\hskip1pt q_i}_{\!\!h_i}$\, and \,$\Gamma^{\hskip1pt q_j}_{\!\!h_j}$\,
bound a cylinder
\,$\mathcal{C}_0\subset{\rm Int}(\mathcal{S})$. In addition, \,$\mathcal{C}_0$\, 
contains a $2$-ball of radius \,$\rho/3$\, since 
\,$d(\Gamma^{\hskip1pt q_i}_{\!\!h_i},\Gamma^{\hskip1pt q_j}_{\!\!h_j})\geq\rho$\, and than
the area \,$area(\mathcal{C}_0)$\, of
\,$\mathcal{C}_0$\, is greater than \,$\delta$.

Now, consider the compact surface \,$\mathcal{S}-{\rm Int}(\mathcal{C}_0)$\, and let 
\,$\mathcal{S}_1\subset\mathcal{S}-{\rm Int}(\mathcal{C}_0)$\, be one of its connected components
whose Euler characteristic is nonzero.
We know that the connected components of 
\,$\partial\mathcal{S}_1$\, are not null homotopic in \,$\Sigma$. Thus, \,$\partial\mathcal{S}_1$\,
has no more that \,$N$\,
connected components. Now, let us choose \,$N$\,
diffeomorphisms \,$\phi_1,\ldots,\phi_N$\, in the list
\,$h_1,\ldots,h_{2N}$\,  which are different from those  used to construct the character curves in
the boundary components of \,$\mathcal{S}_1$. Again, from the induction assumption and by using the
above construction we obtain,  for each \,$i\in\{1,\ldots,N\}$:

\begin{itemize}

\item 
a  point \,$p(\phi_i)\in{\rm Int}(\mathcal{S}_1)$\,  which is a fixed point for all the
diffeomorphisms
\,$f_1,\ldots,f_N,h_1,\ldots,h_{2N+n+1}$\,
except for the diffeomorphism \,$\phi_i$\, and such that \,$p(\phi_i)$\, is a 
\,$\omega$-recurrent point for \,$\phi_i$\,;

\item
a character curve \,$\Gamma^{\hskip1pt p(\phi_i)}_{\!\!\phi_i}\subset{\rm Int}(\mathcal{S}_1)$\, of
\,$\phi_i$\, at  \,$p(\phi_i)$.

\end{itemize}

Furthermore, since  \,$\phi_1,\ldots,\phi_N$\, are in the list
\,$h_1,\ldots,h_{2N}$\, 
it follows  that the distance between any two distinct character curves of the list
$$\Gamma^{\hskip1pt p(\phi_1)}_{\!\!\phi_1},\ldots,\Gamma^{\hskip1pt p(\phi_N)}_{\!\!\phi_N}
\subset\text{Int}(\mathcal{S}_1)$$
is greater than or equal to \,$\rho$. Once more, the topology of \,$\mathcal{S}_1$\, implies that
there exist two distinct curves 
\,$\Gamma^{\hskip1pt p(\phi_i)}_{\!\!\phi_i}$\, and \,$\Gamma^{\hskip1pt p(\phi_j)}_{\!\!\phi_j}$\,
which are homotopic and do not  bound  disks in \,$\mathcal{S}_1$. Thus, they
bound  a cylinder \,$\mathcal{C}_1\subset{\rm Int}(\mathcal{S}_1)$\, 
such that \,$area(\mathcal{C}_1)>\delta$\, since
\,$d(\Gamma^{\hskip1pt p(\phi_i)}_{\!\!\phi_i},
\Gamma^{\hskip1pt p(\phi_j)}_{\!\!\phi_j}\big)\geq\rho$.

Consider now the compact surface \,$\mathcal{S}_1-{\rm Int}(\mathcal{C}_1)$\, and let
\,$\mathcal{S}_2\subset \mathcal{S}_1-{\rm Int}(\mathcal{C}_1)$\, be one of its connected components
whose Euler characteristic is nonzero.  
Once again, we know that 
\,$\partial\mathcal{S}_2$\, has no more than \,$N$\, connected components since they are not null
homotopic in
\,$\Sigma$.  Now, consider
\,$N$\, diffeomorphisms
\,$\psi_1,\ldots,\psi_N$\, in the list  
\,$h_1,\ldots,h_{2N}$\, different from those used to construct the character curves in
the boundary of \,$\mathcal{S}_2$. Repeating the construction above we obtain
for each \,$i\in\{1,\ldots,N\}$: 

\begin{itemize}

\item 
a  point \,$p(\psi_i)\in{\rm Int}(\mathcal{S}_2)$\,  which is a fixed point for all the
diffeomorphisms
\,$f_1,\ldots,f_N\,,h_1,\ldots,h_{2N+n+1}$\,
except for the diffeomorphism \,$\psi_i$\, and such that \,$p(\psi_i)$\, is a 
\,$\omega$-recurrent point for \,$\psi_i$\,;

\item
a character curve \,$\Gamma^{\hskip1pt p(\psi_i)}_{\!\!\psi_i}$\, for \,$\psi_i$\, at 
\,$p(\psi_i)$\, such that, two character curves satisfy
\,$d\big(\Gamma^{\hskip1pt p(\psi_i)}_{\!\!\psi_i}, \Gamma^{\hskip1pt
p(\psi_j)}_{\!\!\psi_j}\big)\geq\rho$\, for all \,$i\neq j$\, and \,$i,j\in\{1,\ldots,N\}$.

\end{itemize}

Applying this construction successively we obtain an infinite family of pairwise disjoint cylinders 
\,$(\mathcal{C}_i)_{i\geq0}$\, in \,${\rm Int}(\mathcal{S})$\, such that the area of each cylinder
is greater than \,$\delta$. This is a contradiction since the area of 
\,$\Sigma$\,  is finite.
\end{proof}

\medskip

Applying Theorem \ref{TPF} for \,$\mathcal{S}=\Sigma$\, we obtain:


\medskip

\begin{thm}\label{TF}

If \,$f_1,\ldots,f_n\in\mathcal{V}_{3N+1}$\, are commuting
diffeomorphisms of \,$\Sigma$\, then
they have common fixed points. 

\end{thm}


\medskip

Now, Theorem \ref{TC} follows from the above result and from the following reasoning 
which can be found in \cite{l2} and \cite{dff}.

Let \,$\mathcal{F}\subset\mathcal{V}_{3N+1}$\, be a nonempty set of commuting diffeomorphisms of
\,$\Sigma$\, and let \,$\emptyset\neq\mathcal{G}\subset\mathcal{F}$\, be a finite subset of 
\,$\mathcal{F}$. Then, by Theorem \ref{TF} it follows that
\,$\text{Fix}(\mathcal{G})\neq\emptyset$. Hence, the family 
\,$\{\text{Fix}(f)\}_{f\in\mathcal{F}}$\,  of closed subsets
of \,$\Sigma$\, has the ``finite intersection property'' which implies that 
\,$\text{Fix}(\mathcal{F})\neq\emptyset$\, and proves Theorem \ref{TC}.


\vskip35pt

\section{The main lemma}
\vskip3mm


Now we prove the Main Lemma which is necessary to obtain Lemma \ref{V.n} according to our strategy.
We shall prove that the \,$k+1$\, diffeomorphisms
$$f_0 \, , \, \ldots \, , \, f_{k-1} \, , \, 
f_k^{\tau_k}\mbox{\tiny$\circ$} \cdots \mbox{\tiny$\circ$} \, f_{3N}^{\tau_{3N}}$$
have a common fixed point in  \,$\text{Int}(\mathcal{S})$\, provided that convenient character curves
in the boundary of \,$\mathcal{S}$\, and convenient exponents 
\,$\tau_k \, , \, \ldots \, ,  \tau_{3N}$\, are available.

For the sake of simplicity we use \,$f_0$\, to denote the identity map of the surface
\,$\Sigma$.

The proof of the Main Lemma will be by induction on \,$k$\, and by contradiction.
We do it in two steps. The first one is 
Lemma \ref{F.1} where we prove the case \,$k=1$, that is, we prove that  
\,$f_1^{\tau_1}\mbox{\tiny$\circ$} \cdots \mbox{\tiny$\circ$} \, f_{3N}^{\tau_{3N}}$\, has 
a fixed point in \,$\text{Int}(\mathcal{S})$\,  by means of the classical Poincar\'e Theorem for
singularities of vector fields.

The second step is the  Main Lemma itself
where we prove the case \,$k>1$. 
It follows very closely the proof of Theorem \ref{TPF}.

In order to simplify the statement of the next lemmas we introduce the following definition.

\medskip

\begin{defn}\label{defchageral}
Let \,$1\leq k\leq 3N-1$\, be an integer. 
Consider subsets \,$\Lambda_k\,,\ldots,\Lambda_{3N}\subsetneq\{1,...,2N\}$\, and integers
\,$\tau_j\in\{1,\ldots,2N\}-\Lambda_j$\, for all \,$j\in\{k,\ldots,3N\}$. Given a subset 
\,$\Lambda\subset\{1,...,3N\}$\, together with commuting diffeomorphisms
\,$f_1,\ldots,f_{3N}$\, one says that a simple  closed curve
\,$\Gamma\subset \Sigma$\, is a \emph{character curve associated to} \,$\Lambda_k
,\ldots,\Lambda_{3N}$\,, $\tau_k,\ldots,\tau_{3N}$\,, $\Lambda$\,, $f_1,\ldots,f_{3N}$\, if\,:
\begin{itemize}
\item [(1)]
either there is \,$\lambda\in\Lambda$\, and a point 
$$p\in \text{Fix}(f_1,\ldots,\widehat{f_{\lambda}},\ldots,f_{3N})-\text{Fix}(f_{\lambda})$$
for which \,$\Gamma$\, is a character curve of \,$f_{\lambda}$\, at \,$p$\,;

\item[(2)]
or there are \,$\xi\in\{k,\ldots,3N-1\}$\,, a number \,$i\in\Lambda_{\xi}$\, 
 and a point 
$$\mu\in\text{Fix}(f_0,\ldots,f_{\xi-1},
f_{\xi}^i \mbox{\tiny$\circ$} f_{\xi+1}^{\tau_{\xi+1}}
\mbox{\tiny$\circ$} \cdots \mbox{\tiny$\circ$} f_{3N}^{\tau_{3N}})
-\text{Fix}(f_{\xi})$$
such that \,$\Gamma$\, is a character curve of \,$f_{\xi}$\, at \,$\mu$.

\end{itemize}

\end{defn}

\medskip

Now, we have the following lemma.

\medskip

\begin{lemma}\label{F.1}

Consider subsets \,$\Lambda_1,\ldots,\Lambda_{3N} \subsetneq \{1,\ldots,2N\}$, integers 
\,$\tau_j\in\{1,\ldots,2N\}-\Lambda_j$\, for all \,$j\in\{1,\ldots,3N\}$\,  and a subset
\,$\Lambda \subset \{1,\ldots,3N\}$\,    such that
$$0\leq \#(\Lambda_1)+\cdots+\#(\Lambda_{3N})+\#(\Lambda)\leq N.$$
Let \,$f_1,\ldots,f_{3N}\in\mathcal{V}_{3N+1}$\, be commuting diffeomorphisms  
and let    \,$\mathcal{S}\subset\Sigma$\, be a compact connected surface with
\,$\chi(\mathcal{S})\neq0$.  Suppose that  any connected component \,$\Gamma$
of 
\,$\partial \mathcal{S}$\, satisfies the following$\,:$ 

\begin{itemize}
\item 
$\Gamma$ is   not null homotopic 
in \,$\Sigma\,;$ 

\item
$\Gamma$  is a character curve associated to
\,$\Lambda_1,\ldots,\Lambda_{3N}$\,, $\tau_1,\ldots,\tau_{3N}$\,, $\Lambda$\,, $f_1,\ldots,f_{3N}$.

\end{itemize}

Then  \ $f_1^{\tau_1} \mbox{\tiny$\circ$} \, f_2^{\tau_2} \mbox{\tiny$\circ$}
\cdots
\mbox{\tiny$\circ$} \, f_{3N}^{\tau_{3N}}$\,  has a fixed point in \,${\rm Int}(\mathcal{S})$.

\end{lemma}


\begin{proof} 

We have that \,$1\leq \tau_1\leq \cdots\leq\tau_{3N}\leq 2N$\, and 
\,$f_1,\ldots,f_{3N}\in  \mathcal{V}_{3N+1}$.
Then, it follows from the definition of the decreasing nested sequence 
\,$\{\mathcal{V}_{k}\}_{k\geq1}$\, that \,$f_i\in\mathcal{V}_{3N+1}\subset\mathcal{V}_{2}$\,
and 
$$f_{1}^{\tau_1} \mbox{\tiny$\circ$} \, \cdots \mbox{\tiny$\circ$} \, 
\widehat{f_{i}^{\tau_i}}  \mbox{\tiny$\circ$} \, \cdots \mbox{\tiny$\circ$} \,
f_{3N}^{\tau_{3N}}\in \mathcal{V}_{3N+1-(3N-1)}=\mathcal{V}_{2} \,.$$ 
Now, from the commutativity of \,$f_1,\ldots,f_{3N}$\, we conclude that
$$f_{1}^{\tau_1}  \mbox{\tiny$\circ$} \, \cdots \mbox{\tiny$\circ$} \,
f_{3N}^{\tau_{3N}}=
\big(f_{1}^{\tau_1} \mbox{\tiny$\circ$} \, \cdots \mbox{\tiny$\circ$} \, 
\widehat{f_{i}^{\tau_i}}  \mbox{\tiny$\circ$} \, \cdots \mbox{\tiny$\circ$} \,
f_{3N}^{\tau_{3N}}\big) \mbox{\tiny$\circ$} \, f_{i}^{\tau_i}$$
has the form \,$h \,\mbox{\tiny$\circ$} \, f^{\ell}$\, 
where  \,$h \,,f \in \mathcal{V}_{2}$\,, 
\,$f\in\{f_1,\ldots,f_{3N}\}$\, and \,$1\leq \ell\leq 2N$.

Let us consider the first type of connected component of \,$\partial\mathcal{S}$, described
in item (1) of Definition \ref{defchageral}.

From item (\ref{char}) of Lemma \ref{bona2} we know that the character curve 
\,$\Gamma^{p}_{\!\!f_\lambda}$\, is \,$\frac{\pi}{10}$-tangent to 
the vector field \,$X_{\!f_{\lambda}}$. Then, it follows from item 
(\ref{fixp2})  of Lemma \ref{bona2}  (case \,$i=0$\, in item (\ref{fixp2}))  
that  the connected components of \,$\partial\mathcal{S}$\, of  type  
\,$\Gamma^{p}_{\!\!f_\lambda}$\, are \,$\frac{2\pi}{5}$-tangent
to the vector field 
\,$X_{\!f_1^{\tau_1}\mbox{\tiny$\circ$} \, \cdots \, \mbox{\tiny$\circ$}
f_{\lambda}^{{\tau}_{\lambda}} \mbox{\tiny$\circ$} \,
\cdots \, \mbox{\tiny$\circ$} f_{3N}^{{\tau}_{3N}}}$.

For the second type of 
connected component of \,$\partial\mathcal{S}$\, we have the following.
If 
\,$\Gamma^{\mu}_{\!\!f_{\xi}}$\, is a connected component  of 
\,$\partial\mathcal{S}$\, for some integer \,$\xi\in\{1,\ldots,3N-1\}$\, then
Definition \ref{defchageral} implies that
$$\mu\in\text{Fix}(
f_1^{\tau_1} \mbox{\tiny$\circ$}\cdots
\mbox{\tiny$\circ$} f_{\xi}^{i} \mbox{\tiny$\circ$} f_{\xi+1}^{\tau_{\xi+1}}
\mbox{\tiny$\circ$} \cdots \mbox{\tiny$\circ$} f_{3N}^{\tau_{3N}}) -\text{Fix}(f_{\xi})$$ where
\,$i\in\Lambda_{\xi}\subset\{1,\ldots,2N\}$. Besides,  item (\ref{char}) of Lemma
\ref{bona2} yields that 
\,$\Gamma^{\mu}_{\!\!f_{\xi}}$\, is \,$\frac{\pi}{10}$-tangent 
to the vector field \,$X_{\!f_{\xi}}$. Therefore, it follows from item 
(\ref{fixp2})  of Lemma \ref{bona2} that  \,$\Gamma^{\mu}_{\!\!f_{\xi}}$\, 
is \,$\frac{2\pi}{5}$-tangent 
to the vector field
\,$X_{\!f_1^{\tau_1} \mbox{\tiny$\circ$} \, \cdots 
\, \mbox{\tiny$\circ$} f_{\xi}^{{\tau}_{\xi}} \mbox{\tiny$\circ$} \, \cdots \, \mbox{\tiny$\circ$}
f_{3N}^{{\tau}_{3N}}}$\, since \,$\tau_{\xi}\in\{1,\ldots,2N\}-\Lambda_{\xi}$.

Thus, 
\,$X_{\!f_1^{\tau_1} \mbox{\tiny$\circ$} \, \cdots \, \mbox{\tiny$\circ$} f_{3N}^{{\tau}_{3N}}}$\,
 has a singularity
in \,$\text{Int}(\mathcal{S})$\, by the classical Poin\-ca\-r\'e Theorem, since
\,$\chi(\mathcal{S})\neq0$.
Hence, the map \,$f_1^{\tau_1} \mbox{\tiny$\circ$} \cdots 
\mbox{\tiny$\circ$} \, f_{3N}^{{\tau}_{3N}}$\, has a
fixed point in \,$\text{Int}(\mathcal{S})$\, and the proof is finished. 
\end{proof}

\medskip

\begin{ml}\label{G.1}

Consider subsets   \,$\Lambda_k,\ldots,\Lambda_{3N} \subsetneq \{1,\ldots,2N\}$,
integers \,$\tau_j\in\{1,\ldots,2N\}-\Lambda_j$\, for all \,$j\in\{k,\ldots,3N\}$\,  and
a subset \,$\Lambda \subset \{1,\ldots,3N\}$\,   such that
$$0\leq \#(\Lambda_k)+\cdots+\#(\Lambda_{3N})+\#(\Lambda)\leq N
\quad \text{where} \quad 1\leq k\leq 3N-1\,.$$
Let \,$f_1,\ldots,f_{3N}\in\mathcal{V}_{3N+1}$\, be commuting diffeomorphisms  
and let    \,$\mathcal{S}\subset\Sigma$\, be a compact connected surface
with \,$\chi(\mathcal{S})\neq0$. 
Suppose that  any connected component \,$\Gamma$ of
\,$\partial \mathcal{S}$\, satisfies the following$\,:$ 

\begin{itemize}
\item 
$\Gamma$ is not null homotopic in \,$\Sigma\,;$

\item 
$\Gamma$ is  a character curve associated to
\,$\Lambda_k\,,\ldots,\Lambda_{3N}$\,, $\tau_k,\ldots,\tau_{3N}$\,, $\Lambda$\,, 
$f_1,\ldots,f_{3N}$.
\end{itemize}

Then the diffeomorphisms \,$f_0 \,, \ldots \,, f_{k-1} \,, f_k^{\tau_k} \mbox{\tiny$\circ$} \cdots
\mbox{\tiny$\circ$} \, f_{3N}^{\tau_{3N}}$\,  have common fixed points in \,${\rm
Int}(\mathcal{S})$.

\end{ml}


\medskip

To prove, for example, that the case \,$k=1$\, implies the case \,$k=2$\, in the Main Lemma, 
we need to construct 
\,$N$\, convenient  character curves supported by some special positive semi-orbits
as in the proof of Theorem \ref{TPF}. For this  we
need to prove again that these semi-orbits stay in \,$\text{Int}(\mathcal{S})$\, and that the
corresponding character curves are pairwise disjoint and do not intersect the character curves in the
boundary of
\,$\mathcal{S}$. These two fundamental steps are proved in the next
two lemmas that we state without proof. Their proofs will be given in the last section. These two
lemmas are actually a blend of Lemma 4.1 of \cite[page~67]{b1} and Lemma 4.2 of
\cite[page~106]{b2} in a more general settings.

\medskip

We recall that the notation \,$f_1,\ldots,\widehat{f_{\lambda}},\ldots, f_{m}$\, means that
\,$f_{\lambda}$\, is not in the  list. Similarly, \,$k_j$\, is not in the list 
\,$k_1,\ldots,\widehat{k_j}, \ldots, k_{n}$.

\medskip

\begin{lemma}\label{fixp4} 
Let \,$f_1,\ldots,f_{3N}\in\mathcal{V}_{3N+1}$\, be
commuting diffeomorphisms  and let 
\,$1\leq \lambda\,,\,\xi\leq 3N$\, be  integers.
Let 
\,$\Gamma^{\hskip1pt p_{\lambda}}_{\!\!f_{\lambda}}$\, and \,$\Gamma^{\hskip1pt
{\mu_j}}_{\!\!f_{\xi}}$\, be  supported by
\,$\mathcal{O}^{+}_{p_{\lambda}}(f_{\lambda})$\, and 
\,$\mathcal{O}^{+}_{\mu_j}(f_{\xi})$\, 
respectively where

\begin{itemize}
\item 
$p_\lambda\in {\rm
Fix}(f_1,\ldots,{\widehat f_\lambda},\ldots,f_{3N})-
{\rm Fix}(f_\lambda)\, ;$ 

\item 
\,$\mu_j\in {\rm Fix}(f_1^{k_1} \mbox{\tiny$\circ$}
 \cdots \mbox{\tiny$\circ$} \, f_{\xi}^{j} \mbox{\tiny$\circ$} \cdots
\mbox{\tiny$\circ$} \, f_{3N}^{k_{3N}})-{\rm
Fix}(f_{\xi}) \, ;$ 
\end{itemize}
and \,$j,k_1,\ldots,\widehat{k_{\xi}},\ldots,k_{3N}\in\{1,\ldots,2N\}$.
Then,  we have\,$:$

\begin{itemize}
\item  
$\gamma^{\hskip1pt \mu_j}_{f_{\xi}} \, \cap \, \Gamma^{\hskip1pt
p_\lambda}_{\!\!f_\lambda}=\emptyset \, ;$

\item
$\Gamma^{\hskip1pt {\mu_j}}_{\!\!f_{\xi}}\cap\Gamma^{\hskip1pt p_\lambda}_{\!\!f_\lambda}=\emptyset
\,;$

\item
$\overline{\mathcal{O}^{+}_{\mu_j}(f_{\xi})} \, \cap \, 
\Gamma^{\hskip1pt p_\lambda}_{\!\!f_\lambda}=\emptyset$.

\end{itemize}

Moreover,
let \,$\rho>0$\, and \,$i\neq j$\, with \,$i\,,j\in\{1,\ldots,2N\}$\,  be such that
$$d\big(x\,,f_{\xi}(x)\big) \geq \rho \ , \ \forall \, x\in 
\mathcal{O}^{+}_{\mu_i}(f_{\xi}) \cup \mathcal{O}^{+}_{\mu_j}(f_{\xi}).$$
Then  \,$d(\Gamma^{\hskip1pt \mu_i}_{\!\!f_{\xi}},\Gamma^{\hskip1pt \mu_j}_{\!\!f_{\xi}})\geq\rho$.

\end{lemma}


\medskip

\begin{lemma}\label{fixp5} 
Let \,$f_1,\ldots,f_m\in\mathcal{V}_{3N+1}$\, be commuting diffeomorphisms
where \,$3\leq m\leq 3N$ 
and let 
\,$\Gamma^{\hskip1pt \nu_j}_{\!\!f_1},\Gamma^{\hskip1pt {\mu_i}}_{\!\!f_{2}}$\, be
  curves supported by 
\,$\mathcal{O}^{+}_{\nu_j}(f_1)$\,, $\mathcal{O}^{+}_{\mu_i}(f_{2})$\, 
respectively satisfying$\,:$

\begin{itemize}
\item 
$\nu_j\in {\rm
Fix}(f_{1}^{j} \mbox{\tiny$\circ$} \, f_2^{k} \mbox{\tiny$\circ$} \,
f_3^{k_{3}} \mbox{\tiny$\circ$} \cdots \mbox{\tiny$\circ$} \, f_{m}^{k_{m}})-
{\rm Fix}(f_1)$\, with \,$1\leq j\leq 2N \, ;$ 

\item 
\,$\mu_i\in {\rm Fix}(f_1\,,f_2^{i} \, \mbox{\tiny$\circ$} \, f_3^{k_3}
\mbox{\tiny$\circ$} \cdots \mbox{\tiny$\circ$} \,  f_{m}^{k_{m}})-{\rm
Fix}(f_{2})$\,  with \,$1\leq i\leq 2N$\, and \,$i\neq k \, ; $ 
\end{itemize}
where \,$k,k_3,\ldots,k_{m}\in\{1,\ldots,2N\}$.
Then,  we have\,$:$

\begin{itemize}

\item
$\gamma^{\nu_j}_{f_1} \, \cap \, \gamma^{\mu_i}_{f_{2}} =\emptyset $
 \ and \
$\Gamma^{\nu_j}_{\!\!f_1} \, \cap \, \Gamma^{\mu_i}_{\!\!f_{2}} =\emptyset \, ;$

\item  
$\gamma^{\nu_j}_{f_1} \, \cap \, \Gamma^{\mu_i}_{\!\!f_{2}} =\emptyset $ \ and \
$\Gamma^{\nu_j}_{\!\!f_1} \, \cap \, \gamma^{\mu_i}_{f_{2}} =\emptyset \, ;$

\item
$\overline{\mathcal{O}^{+}_{\nu_j}(f_1)} \, \cap \, 
\Gamma^{\hskip1pt \mu_i}_{\!\!f_{2}}=\emptyset$ \ and \ 
$\Gamma^{\hskip1pt \nu_j}_{\!\!f_{1}} \cap \,
\overline{\mathcal{O}^{+}_{\mu_i}(f_{2})} = \emptyset$.

\end{itemize}

\end{lemma}


\medskip

The next proof follows very closely the proof of Theorem \ref{TPF}.


\vskip20pt

\subsection{Proof of Main Lemma}


The proof will be by induction on \,$k\in\{1,\ldots,3N-1\}$.
The case \,$k=1$\, is proved in Lemma \ref{F.1}. 

Let us assume that the statement holds for some  \,$k\in\{1,\ldots,3N-2\}$. We will prove that it
also holds  for \,$k+1$.

For this, suppose for a contradiction that the diffeomorphisms 
$$f_1 \,, \, \ldots \,, f_k \quad  \text{and}  \quad 
f_{k+1}^{\tau_{k+1}} \mbox{\tiny$\circ$} \cdots \mbox{\tiny$\circ$} \, f_{3N}^{\tau_{3N}}$$
have no common fixed points in \,$\text{Int}(\mathcal{S})$.
For each \,$j\in\{1,\ldots,2N\}$\, let us consider the maps 
$$f_0 \,, \, \ldots \,, f_{k-1} \quad  \text{and}  \quad 
f_{k}^j \, \mbox{\tiny$\circ$} \, f_{k+1}^{\tau_{k+1}} 
\mbox{\tiny$\circ$} \cdots \mbox{\tiny$\circ$} \, f_{3N}^{\tau_{3N}}.$$

\noindent
The induction assumption on \,$k$\, asserts that these maps have a common fixed point 
 \,$\mu_{k,j}\in\text{Int}(\mathcal{S})$\, for all \,$j\in\{1,\ldots,2N\}$\,
since we can take 
\,$\Lambda_k=\emptyset$\, to apply the induction assumption. Moreover, we have:

\begin{itemize}

\item 
$f_{k}(\mu_{k,j})\neq\mu_{k,j}$\, for all 
\,$j\in\{1,\ldots,2N\}$\, since the maps
\,$f_1,\ldots,f_k$\, and 
\,$f_{k+1}^{\tau_{k+1}} \mbox{\tiny$\circ$} \cdots \mbox{\tiny$\circ$} \,f_{3N}^{{\tau}_{3N}}$\,
have no common fixed points in
\,${\rm Int}(\mathcal{S})$\,;

\item
$\overline{\mathcal{O}^{+}_{\mu_{k,j}}(f_{k})}\subset {\rm Int}(\mathcal{S})$\, 
which follows  from 
   Lemma \ref{fixp4},  and from Lemma   \ref{fixp5} since 
\,$\tau_{\xi}\in\{1,\ldots,2N\}-\Lambda_{\xi}$\, for all integer 
\,$\xi\in\{k+1,\ldots,3N-1\}$.
Lemma \ref{fixp4} guarantees us that \,$\gamma_{\!f_k}^{\hskip1pt \mu_{k,j}}$\, does not intersect
the connected components \,$\Gamma_{\!\!f_{\lambda}}^{p}$\, of \,$\partial\mathcal{S}$, described 
in item (1) of Definition \ref{defchageral}.
Lemma \ref{fixp5} assures that 
\,$\gamma_{\!f_k}^{\hskip1pt \mu_{k,j}}$\, does not intersect
the connected components \,$\Gamma_{\!\!f_{\xi}}^{\mu}$\, 
of \,$\partial\mathcal{S}$, described 
in item (2) of Definition \ref{defchageral}
 for all 
\,$\xi\in\{k+1,\ldots,3N-1\}$.

\end{itemize}

By using Zorn's Lemma and the commutativity of the diffeomorphisms \,$f_1,\ldots,f_{3N}$\,
one can suppose without loss of generality
that
\,$\mu_{k,j}$\, is a
\,$\omega$-recurrent point for \,$f_{k}$. In that case, let 
\,$\Gamma^{\mu_{k,j}}_{\!\!f_{k}}$\, be a character curve for \,$f_{k}$\,  at \,$\mu_{k,j}$. 
By using  Lemmas  \ref{fixp4} and  \ref{fixp5} as above
we conclude that 
\,$\Gamma^{\mu_{k,j}}_{\!\!f_{k}}\subset {\rm Int}(\mathcal{S})$.

By assumption, the maps 
\,$f_1\,,\ldots\,,f_k$\, and  
\,$f_{k+1}^{\tau_{k+1}} \mbox{\tiny$\circ$} \cdots \mbox{\tiny$\circ$} \, f_{3N}^{{\tau}_{3N}}$\,
 have no common fixed points in
   \,${\rm Int}(\mathcal{S})$. Moreover, it follows from item (\ref{fixp}) of Lemma \ref{bona2} 
that the map
$$f_1 \mbox{\tiny$\circ$} \cdots \mbox{\tiny$\circ$} \, f_k \,
\mbox{\tiny$\circ$} \, f_{k+1}^{\tau_{k+1}} \mbox{\tiny$\circ$} \cdots
\mbox{\tiny$\circ$} \, f_{3N}^{{\tau}_{3N}}$$ 
do not have fixed points over \,$\partial\mathcal{S}$.
Thus, 
\,$f_1\,,\ldots\,,f_k$\, and  
\,$f_{k+1}^{\tau_{k+1}} \mbox{\tiny$\circ$} \cdots \mbox{\tiny$\circ$} \, f_{3N}^{{\tau}_{3N}}$\,
 have no common fixed points over \,$\partial\mathcal{S}$.
It results that 
\,$f_1\,,\ldots\,,f_k$\, and  
\,$f_{k+1}^{\tau_{k+1}} \mbox{\tiny$\circ$} \cdots
\mbox{\tiny$\circ$} \, f_{3N}^{{\tau}_{3N}}$\, do not have common
fixed points in \,$\mathcal{S}$.
In such case, there exists \,$\rho>0$\,
satisfying the following condition:
$$d\big(x\,,f_{k}(x)\big) \geq \rho $$
for all 
\,$x\in {\rm
Fix}(f_0\,,\ldots, f_{k-1}\,,f_{k}^{\ell} \, \mbox{\tiny$\circ$} \, f_{k+1}^{\tau_{k+1}}
\mbox{\tiny$\circ$} \cdots \mbox{\tiny$\circ$} \,
f_{3N}^{{\tau}_{3N}})\cap \,\mathcal{S}$\, and for all integer \,$\ell$\, such that
\,$1\leq \ell\leq 2N$.

Let \,$\delta>0$\, be such that the volume of any ball 
of radius  \,$\rho/3$\, contained in \,$\Sigma$\, is greater than \,$\delta$.

Now, from Lemma \ref{fixp4}  we have that the distance between any two
distinct curves of the list of character curves
$$\Gamma^{\hskip1pt \mu_{k,1}}_{\!\!f_{k}},
\ldots,\Gamma^{\hskip1pt \mu_{k,2N}}_{\!\!f_{k}} \subset{\rm Int(\mathcal{S})}$$
is greater than or equal to \,$\rho$. 
Moreover, the topology of \,$\mathcal{S}\subset\Sigma$\, implies that in the  list 
\,$\Gamma^{\hskip1pt \mu_{k,1}}_{\!\!f_{k}},
\ldots,\Gamma^{\hskip1pt \mu_{k,N}}_{\!\!f_{k}}$\, of \,$N$\, elements
 there exist
two distinct curves 
\,$\Gamma^{\hskip1pt \mu_{k,i}}_{\!\!f_{k}}$\, and 
\,$\Gamma^{\hskip1pt \mu_{k,j}}_{\!\!f_{k}}$\,  which are homotopic. On the other hand,  
it follows from Lemma \ref{bondiscsup} that 
these curves can not bound  disks in \,$\mathcal{S}$. 
Thus, we conclude that  these two curves
 bound a cylinder \,$\mathcal{C}_0\subset\text{Int}(\mathcal{S})$\, which contains a ball of radius 
\,$\rho/3$\, since  \,$d(\Gamma^{\hskip1pt \mu_{k,i}}_{\!\!f_{k}},
\Gamma^{\hskip1pt \mu_{k,j}}_{\!\!f_{k}})\geq \rho$. Consequently, 
the area \,$area(\mathcal{C}_0)$\, of 
\,$\mathcal{C}_0$\, is greater than \,$\delta$.

Consider the compact surface \,$\mathcal{S}-\text{Int}(\mathcal{C}_0)$\, and let 
\,$\mathcal{S}_1\subset\mathcal{S}-\text{Int}(\mathcal{C}_0)$\, be one of its connected components
with nonzero Euler characteristic. We know that the connected components of 
\,$\partial\mathcal{S}_1$\, are not null homotopic in \,$\Sigma$. Thus, \,$\mathcal{S}_1$\, has no
more then \,$N$\, connected components in its boundary. 

At this step the connected components of \,$\partial\mathcal{S}_1$\, are character curves
associated to 
$$\Lambda_{k}^{1}\,,\Lambda_{k+1}^{1}\,,\ldots,\Lambda_{3N}^{1}\,,j,\tau_{k+1},
\ldots, \tau_{3N},\Lambda^1,f_1\,,\ldots,f_{3N}$$
where:
\begin{itemize}
\item 
$\Lambda_{k}^1=\big\{i_k \ ; \ \Gamma^{\mu_{k,i_{k}}}_{\!\!f_{k}}\subset\partial\mathcal{S}_1\big\}
\subset\{1,\ldots,N\}$\,;

\item
$\Lambda_{\xi}^1\subset\Lambda_{\xi}$\, for all 
\,$\xi\in\{k+1,\ldots,3N\}$\, and \,$\Lambda^1\subset\Lambda$\,;

\item
$j\in\{1,\ldots,2N\}-\Lambda_{k}^1$\,.

\end{itemize}
 Of course, \,$\Lambda_{k}^1, \ldots,\Lambda^1_{3N},\Lambda^1$\, are such that
$$\#(\Lambda_k^1)+\cdots+\#(\Lambda^{1}_{3N})+\#(\Lambda^{1})\leq N . $$

Now let us go back to the diffeomorphisms
$$f_0 \,, \, \ldots \,, f_{k-1} \quad  \text{and}  \quad 
f_{k}^j \, \mbox{\tiny$\circ$} \, f_{k+1}^{\tau_{k+1}} 
\mbox{\tiny$\circ$} \cdots \mbox{\tiny$\circ$} \, f_{3N}^{\tau_{3N}}$$
 and let
us take \,$j\in\{1,\ldots,2N\}-\Lambda_{k}^{1}$. 
Note that \,$\#\big(\{1,\ldots,2N\}-\Lambda_{k}^{1}\big)\geq N$.
From the induction assumption it results  that they have common fixed points 
\,$\mu_{k,j}\in\text{Int}(\mathcal{S}_1)$\, for all \,$j\in\{1,\ldots,2N\}-\Lambda_{k}^{1}$.
Repeating the same arguments as above we can assume that \,$\mu_{k,j}$\, is \,$\omega$-recurrent
for \,$f_k$\, since \,$\overline{\mathcal{O}^{+}_{\mu_{k,j}}(f_k)}\subset\text{Int}(\mathcal{S}_1)$.
At this point we only need to verify that \,$\gamma^{\hskip1pt \mu_{k,j}}_{\!f_k}$\, do not intersect
\,$\Gamma^{\mu_{k,i_k}}_{\!\!f_k}$\, which is true since \,$j\in\{1,\ldots,2N\}-\Lambda^{1}_{k}$\, 
and \,$i_k\in\Lambda^{1}_{k}$. Following these arguments we conclude that there are
two distinct  curves 
\,$\Gamma^{\hskip1pt \mu_{k,i}}_{\!\!f_{k}}$\, and \,$\Gamma^{\hskip1pt \mu_{k,j}}_{\!\!f_{k}}$\,
in the family of character curves contained in \,$\text{Int}(\mathcal{S}_1)$
$$\big\{\Gamma^{\hskip1pt \mu_{k,j}}_{\!\!f_{k}}\big\}_{j\in\{1,\ldots,2N\}-\Lambda_{k}^{1}}$$
which are not null homotopic in \,$\mathcal{S}_1$\, and
 bound a cylinder \,$\mathcal{C}_1\subset\text{Int}(\mathcal{S}_1)$. Furthermore,
\,$area(\mathcal{C}_1)>\delta$\,  since we have 
\,$d(\Gamma^{\hskip1pt \mu_{k,i}}_{\!\!f_{k}},
\Gamma^{\hskip1pt \mu_{k,j}}_{\!\!f_{k}})\geq \rho$\, for all
\,$i\,,j\in\{1,\ldots,2N\}-\Lambda_{k}^{1}$\, with \,$i\neq j$.

Consider now the compact surface \,$\mathcal{S}_1-\text{Int}(\mathcal{C}_1)$\, and let 
\,$\mathcal{S}_2\subset\mathcal{S}_1-\text{Int}(\mathcal{C}_1)$\, be one of its connected components
whose Euler characteristic is nonzero. Once again, we know that the connected components of 
\,$\partial\mathcal{S}_2$\, are not null homotopic in \,$\Sigma$\, and, consequently, 
\,$\mathcal{S}_2$\, has no more than \,$N$\, connected components in its boundary.

At this step the connected components of \,$\partial\mathcal{S}_2$\, are character curves
associated to 
$$\Lambda_{k}^{2}\,,\Lambda_{k+1}^{2}\,,\ldots,\Lambda_{3N}^{2}\,,j,\tau_{k+1},
\ldots, \tau_{3N},\Lambda^2,f_1\,,\ldots,f_{3N}$$
where:
\begin{itemize}
\item 
$\Lambda_{k}^2=\big\{i_k \ ; \ \Gamma^{\mu_{k,i_{k}}}_{\!\!f_{k}}\subset\partial\mathcal{S}_2\big\}
\subset\{1,\ldots,2N\}$\,;

\item
$\Lambda_{\xi}^2\subset\Lambda_{\xi}^1$\, for all 
\,$\xi\in\{k+1,\ldots,3N\}$\, and \,$\Lambda^2\subset\Lambda^1$\,;

\item
$j\in\{1,\ldots,2N\}-\Lambda_{k}^2$\,.

\end{itemize}
Of course, \,$\Lambda_{k}^2, \ldots,\Lambda^2_{3N},\Lambda^2$\, are such  that
$$\#(\Lambda_k^2)+\cdots+\#(\Lambda^{2}_{3N})+\#(\Lambda^{2})\leq N .$$

Again, let us consider  the diffeomorphisms
$$f_0 \,, \, \ldots \,, f_{k-1} \quad  \text{and}  \quad 
f_{k}^j \, \mbox{\tiny$\circ$} \, f_{k+1}^{\tau_{k+1}}
\mbox{\tiny$\circ$} \cdots \mbox{\tiny$\circ$} \, f_{3N}^{\tau_{3N}}$$
and let
us take \,$j\in\{1,\ldots,2N\}-\Lambda_{k}^{2}$. 
Note that \,$\#\big(\{1,\ldots,2N\}-\Lambda_{k}^{2}\big)\geq N$.
Once more, the induction assumption on \,$k$\, asserts  that
they have common fixed points 
\,$\mu_{k,j}\in\text{Int}(\mathcal{S}_2)$\, for all \,$j\in\{1,\ldots,2N\}-\Lambda_{k}^{2}$.
Repeating exactly the same arguments as above we can assume that \,$\mu_{k,j}$\, 
is \,$\omega$-recurrent for \,$f_k$\, and we 
obtain two distinct  curves 
\,$\Gamma^{\hskip1pt \mu_{k,i}}_{\!\!f_{k}} , 
\Gamma^{\hskip1pt \mu_{k,j}}_{\!\!f_{k}}$\,
in the family of character curves contained in \,$\text{Int}(\mathcal{S}_2)$
$$\big\{\Gamma^{\hskip1pt \mu_{k,j}}_{\!\!f_{k}}\big\}_{j\in\{1,\ldots,2N\}-\Lambda_{k}^{2}}$$
which are not null homotopic in \,$\mathcal{S}_2$\, and
 bound a cylinder \,$\mathcal{C}_2\subset\text{Int}(\mathcal{S}_2)$. We have also that
\,$area(\mathcal{C}_2)>\delta$\,  since the distance 
\,$d(\Gamma^{\hskip1pt \mu_{k,i}}_{\!\!f_{k}},
\Gamma^{\hskip1pt \mu_{k,j}}_{\!\!f_{k}})\geq \rho$\, for all
\,$i\,,j\in\{1,\ldots,2N\}-\Lambda_{k}^{2}$\, with \,$i\neq j$.

By successively repeating the above construction we obtain a family of pairwise disjoint cylinders
\,$\{\mathcal{C}_i\}_{i\geq0}$\, contained in \,$\Sigma$\, and such that
 \,$area(\mathcal{C}_i)>\delta$\, for all integer \,$i\geq0$.
This is however impossible and therefore completes the proof of the statement.


\vskip35pt

\section{Proof of Lemma \ref{V.n}}\label{proofV.n}
\vskip3mm


Suppose for a contradiction that the diffeomorphisms
\,$f_1 , \ldots , f_{3N}$\, have no common fixed points in \,$\text{Int}(\mathcal{S})$. For each
\,$j\in\{1,\ldots,2N\}$\, let us consider the \,$3N-1$\, diffeomorphisms
$$f_1 \,, \ldots \,, f_{3N-2}\,,f_{3N-1}^j \, \mbox{\tiny$\circ$} \, f_{3N}.$$
Taking \,$k=3N-1$\, and \,$\Lambda_{3N-1}=\emptyset$\, in the Main Lemma we conclude that they have
a common fixed point \,$\mu_{3N-1,j}\in\text{Int}(\mathcal{S})$\,  for all 
integer  \,$j\in\{1,\ldots,2N\}$.
Moreover, we have:

\begin{itemize}

\item 
$f_{3N-1}(\mu_{3N-1,j})\neq\mu_{3N-1,j}$\, for all 
\,$j\in\{1,\ldots,2N\}$\, since the diffeomorphisms
\,$f_1,\ldots,f_{3N}$\,
have no common fixed points in
\,${\rm Int}(\mathcal{S})$\,;

\item
$\overline{\mathcal{O}^{+}_{\mu_{3N-1,j}}(f_{3N-1})}\subset {\rm Int}(\mathcal{S})$\, 
which follows  from    Lemma \ref{fixp4}.

\end{itemize}

Zorn's Lemma and the commutativity of the maps \,$f_1,\ldots,f_{3N}$\,
allow us to suppose without loss of generality
that
\,$\mu_{3N-1,j}$\, is a
\,$\omega$-recurrent point for \,$f_{3N-1}$. In that case, let 
\,$\Gamma^{\mu_{3N-1,j}}_{\!\!f_{3N-1}}$\, be a character curve for \,$f_{3N-1}$\,  at
\,$\mu_{3N-1,j}$.  Once more, from  Lemma  \ref{fixp4} 
we have that 
\,$\Gamma^{\mu_{3N-1,j}}_{\!\!f_{3N-1}}\subset {\rm Int}(\mathcal{S})$.

From now on the proof of the lemma is concluded by repeating the
proof of the Main Lemma: it suffices to substitute \,``$k$''\, by \,``$3N-1$''\, and
the  ``induction
assumption'' by the ``Main Lemma''.


\vskip35pt

\section{Proofs of technical lemmas}\label{prooftechlem}
\vskip3mm


In this section we prove the last two technical lemmas. 
The proofs are similar but more technical than the proof of Lemma \ref{fixp4A} because we treat 
the dynamics of curves supported by positive semi-orbits in a more general settings.


\vskip20pt
\subsection{Proof of Lemma \ref{fixp4}}


Firstly note  that for each  \,$j\in\{1,\ldots,2N\}$\,  the map
\,$f_1^{k_1} \mbox{\tiny$\circ$} \cdots \mbox{\tiny$\circ$} \, f_{\xi}^{j}
 \mbox{\tiny$\circ$}  \cdots  \mbox{\tiny$\circ$} \, f_{3N}^{k_{3N}}$\,
has the form
\,$h \, \mbox{\tiny$\circ$}  f^{\ell}$\, where  the maps
\,$h\,,f\in\mathcal{V}_2$\,, \,$f\in\{f_1,\ldots,f_{3N}\}$\, and the integer
  \,$\ell\in\{1,\ldots,2N\}$. 

Hence,
it follows from item (\ref{fixp}) of Lemma \ref{bona2} (case \,$i=0$\, in item (\ref{fixp})) that 
the diffeomorphism 
\,$f_1^{k_1} \mbox{\tiny$\circ$} \cdots
 \mbox{\tiny$\circ$} \, f_{\xi}^{j} \mbox{\tiny$\circ$}  \cdots
 \mbox{\tiny$\circ$} \, f_{3N}^{k_{3N}}$\, has no fixed
points in the ball
 \,$\textrm{B}\big(f^{m}_{\lambda}(p_{\lambda})\,;
3\,d(f^{m}_{\lambda}(p_{\lambda})\,,f^{m+1}_{\lambda}(p_{\lambda}))\big)$.

Now, let us prove that 
$\Gamma^{\hskip1pt {\mu_j}}_{\!\!f_{\xi}}\cap
\Gamma^{\hskip1pt p_\lambda}_{\!\!f_\lambda}=\emptyset$. 
For this, suppose for a contradiction that
there exist \,$m,n,k,l\geq0$\, such that
$$[f^{m}_{\lambda}(p_{\lambda})\,,f^{n}_{\lambda}(p_{\lambda})]\cap 
[f^{k}_{\xi}(\mu_j)\,,f^{l}_{\xi}(\mu_j)]\neq\emptyset$$
where 
\begin{equation}\label{sec4:lem4.2:sitde2}
\begin{split}
d\big(f^{m}_{\lambda}(p_{\lambda})\,,f^{n}_{\lambda}(p_{\lambda})\big) & \leq 
\mbox{$\frac{3}{2}$} \,
d\big(f^{m}_{\lambda}(p_{\lambda})\,,f^{m+1}_{\lambda}(p_{\lambda})\big) 
\\  d\big(f^{k}_{\xi}(\mu_j)\,,f^{l}_{\xi}(\mu_j)\big) & \leq \mbox{$\frac{3}{2}$} \,
d\big(f^{k}_{\xi}(\mu_j)\,,f^{k+1}_{\xi}(\mu_j)\big).
\end{split}
\end{equation}

\noindent
By the triangle inequality we have  
\begin{align}
 d\big(f^{m}_{\lambda}(p_{\lambda})\,,f^{k}_{\xi}(\mu_j)\big) & \leq   
d\big(f^{m}_{\lambda}(p_{\lambda})\,,f^{n}_{\lambda}(p_{\lambda})\big)+
d\big(f^{k}_{\xi}(\mu_j)\,,f^{l}_{\xi}(\mu_j)\big)
\nonumber \\
 & \leq   
\mbox{$\frac{3}{2}$} \,
d\big(f^{m}_{\lambda}(p_{\lambda})\,,f^{m+1}_{\lambda}(p_{\lambda})\big)+\mbox{$\frac{3}{2}$}
\, d\big(f^{k}_{\xi}(\mu_j)\,,f^{k+1}_{\xi}(\mu_j)\big) 
\nonumber\\
& \leq 3\max\Big\{d\big(f^{m}_{\lambda}(p_{\lambda})\,,f^{m+1}_{\lambda}(p_{\lambda})\big) \,,\,
d\big(f^{k}_{\xi}(\mu_j)\,,f^{k+1}_{\xi}(\mu_j)\big)\Big\}. 
\nonumber
\end{align}

\noindent
Consequently, we have two possibilities:

\begin{itemize}
\item [$\bullet$]
either \,$f^{k}_{\xi}(\mu_j)$\, is  in the   ball 
 \,$\textrm{B}\big(f^{m}_{\lambda}(p_{\lambda})\,;
3\,d(f^{m}_{\lambda}(p_{\lambda})\,,f^{m+1}_{\lambda}(p_{\lambda}))\big)$\, which is however
impossible  as follows from the remarks below: 

\begin{itemize}
\item [--]
$f^{k}_{\xi}(\mu_j)\in {\rm Fix}(f_1^{k_1}
 \mbox{\tiny$\circ$} \cdots \mbox{\tiny$\circ$} \, f_{\xi}^{j}
 \mbox{\tiny$\circ$}  \cdots \mbox{\tiny$\circ$} \, f_{3N}^{k_{3N}})$\, by
commutativity;

\item[--]
the diffeomorphism \,$f_1^{k_1} \mbox{\tiny$\circ$} \cdots
 \mbox{\tiny$\circ$} \, f_{\xi}^{j} \mbox{\tiny$\circ$}  \cdots
 \mbox{\tiny$\circ$} \, f_{3N}^{k_{3N}}$\, has
no fixed points in the ball
 \,$\textrm{B}\big(f^{m}_{\lambda}(p_{\lambda})\,,
3\,d\big(f^{m}_{\lambda}(p_{\lambda})\,,f^{m+1}_{\lambda}(p_{\lambda})\big)$\,;

\end{itemize}

\item[$\bullet$]
or \,$f^{m}_{\lambda}(p_{\lambda})$\, is in the ball
\,$\textrm{B}\big(f^{k}_{\xi}(\mu_j)\,;3\,d(f^{k}_{\xi}(\mu_j)\,,f^{k+1}_{\xi}(\mu_j))\big)$:

\begin{itemize}

\item[--]
If \,$\lambda\neq \xi$\, then we have that \,$f^{m}_{\lambda}(p_{\lambda})$\, is a fixed point 
for \,$f_{\xi}$\, which is a contradiction since  \,$f_{\xi}$\, does not have fixed points in the ball 
\,$\textrm{B}\big(f^{k}_{\xi}(\mu_j)\,;3\,d(f^{k}_{\xi}(\mu_j)\,,f^{k+1}_{\xi}(\mu_j))\big)$\,
where \,$f^{k}_{\xi}(\mu_j)\notin {\rm Fix}(f_{\xi})$\,;

\item[--]
If \,$\lambda=\xi$\, then we have
\,$f_{\lambda}^m(p_{\lambda})\in
\textrm{B}\big(f^{k}_{\lambda}(\mu_j)\,;3\,d(f^{k}_{\lambda}
(\mu_j)\,,f^{k+1}_{\lambda}(\mu_j))\big)$\, and
\,$f_{\lambda}^m(p_{\lambda})\in\text{Fix}(f_1^{k_1}
 \mbox{\tiny$\circ$} \cdots \mbox{\tiny$\circ$} \, {\widehat f_{\lambda}^j}
 \mbox{\tiny$\circ$} \cdots \mbox{\tiny$\circ$} \, f_{3N}^{k_{3N}})$. 
On the other hand we have also that
\,$f^{k}_{\lambda}(\mu_j)\in\text{Fix}
(f_1^{k_1} \mbox{\tiny$\circ$} \cdots  \mbox{\tiny$\circ$} \, {f_{\lambda}^j}
 \mbox{\tiny$\circ$} \cdots \mbox{\tiny$\circ$} \, f_{3N}^{k_{3N}})-\text{Fix}(f_{\lambda})$.
Thus, it follows from item (\ref{fixp}) of Lemma \ref{bona2} that the diffeomorphism
\,$f_1^{k_1} \mbox{\tiny$\circ$} \cdots \mbox{\tiny$\circ$} \, {\widehat f_{\lambda}^j} \,
 \mbox{\tiny$\circ$} \cdots \mbox{\tiny$\circ$} \, f_{3N}^{k_{3N}}$\, has no fixed points in the ball
\,$\textrm{B}\big(f^{k}_{\lambda}(\mu_j)\,;3\,d(f^{k}_{\lambda}
(\mu_j)\,,f^{k+1}_{\lambda}(\mu_j))\big)$\, which is a contradiction.

\end{itemize}

\end{itemize}

\noindent
The preceding discussion has shown that 
\,$\Gamma^{\hskip1pt {\mu_j}}_{\!\!f_{\xi}}\cap\Gamma^{\hskip1pt
p_\lambda}_{\!\!f_\lambda}=\emptyset$. 
Similar arguments prove 
that  \,$\gamma^{\hskip1pt \mu_j}_{f_{\xi}} \, \cap \, \Gamma^{\hskip1pt
p_\lambda}_{\!\!f_\lambda}=\emptyset$
\  for all 
\,$j\in\{1,\ldots,2N\}$.

In addition, 
\,$\overline{\mathcal{O}^{+}_{\mu_j}(f_{\xi})}\subset  
{\rm Fix}(f_1^{k_1} \mbox{\tiny$\circ$} \cdots \mbox{\tiny$\circ$} \, f_{\xi}^j
 \mbox{\tiny$\circ$}  \cdots \mbox{\tiny$\circ$} \,
f_{3N}^{k_{3N}})$\, since the diffeomorphisms commute.  Thus, it follows from 
item (\ref{fixp}) of Lemma \ref{bona2} that \,$\overline{\mathcal{O}^{+}_{\mu_j}(f_{\xi})} \,
\cap\, 
\Gamma^{\hskip1pt p_\lambda}_{\!\!f_\lambda}=\emptyset$\, for all \,$j\in\{1,\ldots,2N\}$\, 
as we have seen in the beginning of the proof. This proves the first part of the
lemma. 

To prove the second part  let \,$\rho>0$\, and \,$i\neq j$\, with
\,$i\,,j\in\{1,\ldots,2N\}$\,  be such that
$$d(x\,,f_{\xi}(x)) \geq \rho \ \ , \ \ \forall \, x\in 
\mathcal{O}^{+}_{\mu_i}(f_{\xi}) \cup \mathcal{O}^{+}_{\mu_j}(f_{\xi})\,.
$$
Suppose  that
\,$d(\Gamma^{\hskip1pt \mu_i}_{\!\!f_{\xi}},\Gamma^{\hskip1pt \mu_j}_{\!\!f_{\xi}})<\rho$. In that
case there exist
\,$m,n,k,l\geq0$\, satisfying (\ref{sec4:lem4.2:sitde2}) and 
two points  
\,$a\in[f^m_{\xi}(\mu_i)\,,f^{n}_{\xi}(\mu_i)]  \ \ \mbox{and} 
\ \  b\in[f^k_{\xi}(\mu_j)\,,f^{l}_{\xi}(\mu_j)] $\, 
such that \,$d(a\,,b)< \rho$.
Therefore,
\begin{align}
d\big(f^m_{\xi}(\mu_i)\,,f^k_{\xi}(\mu_j)\big) & \leq
d\big(f^m_{\xi}(\mu_i)\,,f^{n}_{\xi}(\mu_i)\big)+
d(a\,,b)+d\big(f^k_{\xi}(\mu_j)\,,f^{l}_{\xi}(\mu_j)\big) 
\nonumber\\
 & \leq
\mbox{$\frac{3}{2}$}\,d\big(f^m_{\xi}(\mu_i)\,,f^{m+1}_{\xi}(\mu_i)\big)+\rho+\mbox{$\frac{3}{2}$}\,
d\big(f^k_{\xi}(\mu_j)\,,f^{k+1}_{\xi}(\mu_j)\big) 
\nonumber\\
  & \leq  
4\max\Big\{d\big(f^m_{\xi}(\mu_i)\,,f^{m+1}_{\xi}(\mu_i)\big) \,,\,
d\big(f^k_{\xi}(\mu_j)\,,f^{k+1}_{\xi}(\mu_j)\big)\Big\}. 
\nonumber
\end{align}

\noindent
Consequently,
\begin{itemize}
\item [--]
either \,$f^k_{\xi}(\mu_j)$\, is  in the   ball 
 \,$\textrm{B}\big(f^m_{\xi}(\mu_i)\,,4\,d\big(f^m_{\xi}(\mu_i)\,,f^{m+1}_{\xi}(\mu_i))\big)$\, which
is impossible by item (\ref{fixp}) of Lemma \ref{bona2}  since        
\,$f^k_{\xi}(\mu_j)$\, is a fixed point of 
\,$f_1^{k_1} \mbox{\tiny$\circ$} \cdots \mbox{\tiny$\circ$} \, f_{\xi}^j
 \mbox{\tiny$\circ$}  \cdots \mbox{\tiny$\circ$} \, f_{3N}^{k_{3N}}$\,
and
\,$i\neq j$\,;

\item[--]
or \,$f^m_{\xi}(\mu_i)$\, is in the ball
\,$\textrm{B}\big(f^k_{\xi}(\mu_j)\,,4\,d\big(f^k_{\xi}(\mu_j)\,,f^{k+1}_{\xi}(\mu_j))\big)$\, 
which is impossible by  the same reason.
\end{itemize}
The proof of Lemma \ref{fixp4} is finished.


\vskip20pt
\subsection{Proof of Lemma \ref{fixp5}}


As  remarked in the last proof, for each \,$j\in\{1,\ldots,2N\}$\, the map 
\,$f_1^{j} \, \mbox{\tiny$\circ$} \, f_2^{k} \, \mbox{\tiny$\circ$} \, f_3^{k_{3}}
 \mbox{\tiny$\circ$}  \cdots \mbox{\tiny$\circ$} \, f_{m}^{k_{m}}$\, has the form
\,$h \, \mbox{\tiny$\circ$} f^{\ell}$\, where  
\,$h\,,f\in\mathcal{V}_2$\,, \,$f\in\{f_1,\ldots,f_{m}\}$\, and the integer 
\,$\ell$\, lies in \,$\{1,\ldots,2N\}$.

Suppose that there exist \,$m,n,k,l\geq0$\, such that
$$[f^{m}_{1}(\nu_j)\,,f^{n}_{1}(\nu_j)]\cap 
[f^{k}_{2}(\mu_i)\,,f^{l}_{2}(\mu_i)]\neq\emptyset$$
where 
\begin{align}
d\big(f^{m}_{1}(\nu_j)\,,f^{n}_{1}(\nu_j)\big) & \leq \mbox{$\frac{3}{2}$} \,
d\big(f^{m}_{1}(\nu_j)\,,f^{m+1}_{1}(\nu_j)\big) \nonumber
\\  d\big(f^{k}_{2}(\mu_i)\,,f^{l}_{2}(\mu_i)\big) & \leq \mbox{$\frac{3}{2}$} \,
d\big(f^{k}_{2}(\mu_i)\,,f^{k+1}_{2}(\mu_i)\big).
\nonumber
\end{align}

\noindent
By the triangle inequality we have  
\begin{align}
 d\big(f^{m}_{1}(\nu_j)\,,f^{k}_{2}(\mu_i)\big) & \leq   
d\big(f^{m}_{1}(\nu_j)\,,f^{n}_{1}(\nu_j)\big)+
d\big(f^{k}_{2}(\mu_i)\,,f^{l}_{2}(\mu_i)\big)
\nonumber \\
 & \leq   
\mbox{$\frac{3}{2}$} \, d\big(f^{m}_{1}(\nu_j)\,,f^{m+1}_{1}(\nu_j)\big)+\mbox{$\frac{3}{2}$}
\, d\big(f^{k}_{2}(\mu_i)\,,f^{k+1}_{2}(\mu_i)\big) 
\nonumber\\
& \leq 3\max\Big\{d\big(f^{m}_{1}(\nu_j)\,,f^{m+1}_{1}(\nu_j)\big) \,,\,
d\big(f^{k}_{2}(\mu_i)\,,f^{k+1}_{2}(\mu_i)\big)\Big\}. 
\nonumber
\end{align}

\noindent
Consequently,
\begin{itemize}
\item [$\bullet$]
either \,$f^{k}_{2}(\mu_i)$\, is  in the   ball 
 \,$\textrm{B}\big(f^{m}_{1}(\nu_j)\,;
3\,d\big(f^{m}_{1}(\nu_j)\,,f^{m+1}_{1}(\nu_j))\big)$\, which is
impossible   since \,$f^{k}_{2}(\mu_i)\in{\rm Fix}(f_1)$\, and \,$f_1$\, 
has no fixed points in the ball  
 \,$\textrm{B}\big(f^{m}_{1}(\nu_j)\,;
3\,d\big(f^{m}_{1}(\nu_j)\,,f^{m+1}_{1}(\nu_j)\big)$\, because 
\,$f_1^m(\nu_j)\notin{\rm Fix}(f_1)$\,;

\item[$\bullet$]
or \,$f^{m}_{1}(\nu_j)\in
\textrm{B}\big(f^{k}_{2}(\mu_i)\,;3\,d\big(f^{k}_{2}(\mu_i)\,,f^{k+1}_{2}(\mu_i))\big)$\, 
which is also impossible as follows from the following remarks:

\begin{itemize}
\item 
$f^{m}_{1}(\nu_j)$\, is a fixed point of 
\,$f_1^{j}  \mbox{\tiny$\circ$} f_2^k  \mbox{\tiny$\circ$} f_3^{k_3}
 \mbox{\tiny$\circ$} \cdots \mbox{\tiny$\circ$}  f_{m}^{k_{m}}$\,; 

\item
 the diffeomorphism  
\,$f_1^{j} \mbox{\tiny$\circ$} f_2^k \mbox{\tiny$\circ$}  f_3^{k_3}
 \mbox{\tiny$\circ$} \cdots \mbox{\tiny$\circ$}  f_{m}^{k_{m}}$\, 
has no fixed points in the ball
\,$\textrm{B}\big(f^{k}_{2}(\mu_i)\,;3\,d\big(f^{k}_{2}(\mu_i)\,,f^{k+1}_{2}(\mu_i)\big)$\, 
when \,$i\neq k$\, as a consequence of item (\ref{fixp}) of Lemma \ref{bona2}. 
Note that we have \,$\mu_i\in\text{Fix}(f_1^j
 \mbox{\tiny$\circ$} f_2^i \, \mbox{\tiny$\circ$} f_3^{k_3} \mbox{\tiny$\circ$} \cdots
 \mbox{\tiny$\circ$}  f_m^{k_m})$\,
and \,$\mu_i\notin\text{Fix}(f_2)$.

\end{itemize}

\end{itemize}

\noindent
This argument proves the first two items of the lemma.

To finish the proof we remark that 
\,$\overline{\mathcal{O}^{+}_{\nu_j}(f_1)}\subset{\rm Fix}
(f_1^{j} \mbox{\tiny$\circ$} f_2^k \mbox{\tiny$\circ$}  f_3^{k_3}
 \mbox{\tiny$\circ$} \cdots \mbox{\tiny$\circ$}  f_{m}^{k_{m}})$.
On the other hand, 
\,$\mu_i\in{\rm Fix}(f_1^{j} \, \mbox{\tiny$\circ$} \,f_2^i
\, \mbox{\tiny$\circ$} \, f_3^{k_3} \mbox{\tiny$\circ$} \cdots
 \mbox{\tiny$\circ$} \, f_{m}^{k_{m}})-{\rm
Fix}(f_{2})$\, and consequently, there is no fixed point of 
\,$f_1^{j} \, \mbox{\tiny$\circ$} \, f_2^k \, \mbox{\tiny$\circ$} \, f_3^{k_3}
 \mbox{\tiny$\circ$} \cdots  \mbox{\tiny$\circ$} \,
f_{m}^{k_{m}}$\, over the curve 
\,$\Gamma^{\hskip1pt \mu_i}_{\!\!f_{2}}$\, since \,$i\neq k$. Thus,
\ $\overline{\mathcal{O}^{+}_{\nu_j}(f_1)} \, \cap \, 
\Gamma^{\hskip1pt \mu_i}_{\!\!f_{2}}=\emptyset$.

Finally, 
\ $\Gamma^{\hskip1pt \nu_j}_{\!\!f_{1}} \cap \,
\overline{\mathcal{O}^{+}_{\mu_i}(f_{2})} = \emptyset$ \ since 
\,$\overline{\mathcal{O}^{+}_{\mu_i}(f_{2})}\subset{\rm Fix}(f_1)$\, and \,$f_1$\, 
has no fixed points over \,$\Gamma^{\hskip1pt \nu_j}_{\!\!f_{1}}$. The proof is finished.


\vskip20pt
\subsubsection*{Acknowledgement} It is a pleasure to thank the referees for their comments and
suggestions that allowed me to improve on the first version of this article.
I am especially indebted to the second referee who suggested me Definition 6.1 along with
precise modifications on some lemmas that made this note more readable.



\end{document}